\documentclass[10pt,journal, onecolumn]{IEEEtran}

\usepackage{color,latexsym,amsfonts,amssymb,bbm,comment,amsmath,cite, amsthm}
\usepackage{hyperref}
\usepackage{graphicx}
\usepackage{tikz}
\usetikzlibrary{decorations.pathreplacing}
\usetikzlibrary{patterns}
\usepackage{mathtools}
\usepackage{ifthen}
\usepackage{babel}

\usepackage{adjustbox}
\usepackage{xcolor}
\usepackage{enumerate}
\usepackage{amsmath}
\usepackage{amssymb}
\usepackage{fancyvrb}
\usepackage{grffile}
\usepackage{hyperref}
\usepackage{longtable}
\usepackage{booktabs}
\usepackage[inline]{enumitem}
\usepackage[normalem]{ulem}
\usepackage{float}
\usetikzlibrary{calc} 
\usetikzlibrary{arrows}
\pgfmathsetseed{14}

\newcommand{\coxSegment}[6]{

	\pgfmathsetmacro{\Xa}{#1}
	\pgfmathsetmacro{\Xb}{#2}
	\pgfmathsetmacro{\Xc}{#3}
	\pgfmathsetmacro{\Xd}{#4}

	\draw (\Xa - \factor*\Xc + \factor*\Xa,\Xb  - \factor*\Xd + \factor*\Xb) -- (\Xc + \factor*\Xc - \factor*\Xa, \Xd + \factor*\Xd - \factor*\Xb);
	\foreach \y in {1,2,...,#5}{
		\pgfmathsetmacro{\Xf}{random()}
		\pgfmathsetmacro{\Xg}{1- \Xf)}
		\pgfmathsetmacro{\Xh}{\Xa*\Xf + \Xc*\Xg}
		\pgfmathsetmacro{\Xi}{\Xb*\Xf+ \Xd*\Xg}
		\pgfmathparse{\Xh*\Xh + \Xi*\Xi > #6 ? 1: 0}
		\ifthenelse{\pgfmathresult>0}{\fill (\Xh, \Xi) circle (2pt);}{}
	}

}

\def\eq{\begin{equation}}
\def\en{\end{equation}}

\newtheorem{theorem}{Theorem}[section]

\newtheorem{lemma}[theorem]{Lemma}
\newtheorem{proposition}[theorem]{Proposition}
\theoremstyle{definition}
\newtheorem{definition}[theorem]{Definition}

\newcommand{\muv}{\mu}

\newcommand{\lac}{\la_{\mathrm c}}

\newcommand{\supp}{{\mathrm supp}}
\usepackage{color}

\newcommand{\vmin}{v_{\mathrm{min}}}
\newcommand{\vminp}{v'_{\mathrm{min}}}
\newcommand{\vmax}{v_{\mathrm{max}}}
\newcommand{\la}{\lambda}
\newcommand{\dist}{\mathrm{dist}}
\newcommand{\one}{\mathbbmss{1}}
\def \k{\kappa}

\def\e{{\varepsilon}}
\def\eps{\varepsilon}

\newcommand{\R}{\mathbb R}
\newcommand{\Q}{\mathbb Q}
\newcommand{\E}{\mathbb E}

\def\P{\mathbb P} 
\newcommand{\Z}{\mathbb Z}
\newcommand{\N}{\mathbb N}
\newcommand{\C}{\mathcal{C}}
\def\a{\alpha}
\def\d{{\mathrm d}}

\def\L{\Lambda}
\def\supp{{\mathrm{supp}}}
\newcommand{\laW}{{\lambda_{\rm W}}}
\newcommand{\rhoW}{{\rho_{\rm W}}}
\newcommand{\rhoI}{{\rho_{\rm I}}}

\newcommand{\varrhoW}{{\varrho_{\rm W}}}
\newcommand{\varrhoI}{{\varrho_{\rm I}}}
\newcommand{\varrhoIb}{{\varrho^b_{\rm I}}}

\newcommand{\varrhoImin}{{\varrho_{\rm I,min}}}
\newcommand{\varrhoIminb}{{\varrho^b_{\rm I,min}}}
\newcommand{\varrhoWmin}{{\varrho_{\rm W,min}}}

\newcommand{\laWc}{{\lambda_{\rm W,c}}}

\begin{document}
\title{Chase-escape in dynamic device-to-device networks}	

\author{
Elie~Cali,~\IEEEmembership{Orange,}
Alexander~Hinsen,~\IEEEmembership{WIAS,}
        Benedikt~Jahnel,~\IEEEmembership{WIAS \& TU Braunschweig,}
        and Jean-Philippe~Wary,~\IEEEmembership{Orange}
\thanks{Orange Labs, 44 Avenue de la République, 92320 Châtillon, France}
\thanks{Weierstrass Institute for Applied Analysis and Stochastics, Mohrenstraße 39, 10117 Berlin, Germany}
\thanks{Technische Universit\"at Braunschweig, Institut f\"ur Mathematische Stochastik, Universit\"atsplatz 2, 38106 Braunschweig, Germany}}

\maketitle
		
\begin{abstract}
The present paper features results on global survival and extinction of an infection in a multi-layer network of mobile agents. Expanding on a model first presented in~\cite{HJCW22}, we consider an urban environment, represented by line-segments in the plane, in which agents move according to a random waypoint model based on a Poisson point process. Whenever two agents are at sufficiently close proximity for a sufficiently long time the infection can be transmitted and then propagates into the system according to the same rule starting from a typical device. Inspired by wireless network architectures, the network is additionally equipped with a second class of agents that is able to transmit a patch to neighboring infected agents that in turn can further distribute the patch, leading to a chase-escape dynamics. 

We give conditions for parameter configurations that guarantee existence and absence of global survival as well as an in-and-out of the survival regime, depending on the speed of the devices. We also provide complementary results for the setting in which the chase-escape dynamics is defined as an independent process on the connectivity graph. 
The proofs mainly rest on percolation arguments via discretization and multiscale analysis. 

\end{abstract}




\section{Introduction and setting}\label{sec_intro}
In the past decades we have seen a tremendous increase in demand for data exchange on a global scale creating significant pressure on network operators to maintain a good quality of service. One particularly interesting concept in this development is {\em device-to-device (D2D) communications}. The idea is that network components can exchange data in a peer-to-peer fashion over a wireless channel and thereby constitute a decentralized ad-hoc network. Among the many potential benefits of these networks, such as robustness, communication speed, cost efficiency etc., one of the key challenges is their high complexity and unpredictability. 

In order to cope with these challenges, since the early sixties, probabilistic modeling and analysis has been developed and employed to provide qualitative and quantitative insights for D2D networks, see for example~\cite{haenggi2012stochastic,franceschetti2007random,baccelli2009stochastic1} and many more. In this context, one of the most fruitful approaches is based on {\em stochastic geometry}, where network components are conceived as point point processes in space~\cite{blaszczyszyn2018stochastic,jahnel2020probabilistic}. Such a random point cloud is then often seen as the vertex set of a random spatial graph in which edges are drawn according to some rule that reflects the peer-to-peer communication paradigm. In its purest form this is the classical {\em Poisson--Gilbert graph}~\cite{gilbert61}, which serves as a prototypical model for the spatial clustering behavior of a D2D network in a non-dynamical setting. In particular, already in this model, one can observe the celebrated {\em phase transition of percolation} in the following sense. Let $\lambda>0$ be the intensity of the underlying homogeneous Poisson point process $X$ and $r>0$ the connectivity threshold, i.e., the Poisson--Gilbert graph has vertex set $X$ and any pair of vertices is connected by an edge iff their Euclidean distance is less than $r$. Then, there exists a finite and positive $\la_{\rm c}=\la_{\rm c}(r)$ such that for $\la< \la_{\rm c}$ the graph does not contain any infinite connected component almost surely and for $\la> \la_{\rm c}$ such an infinite connected component exists with probability one. In the context of D2D networks, the parameter regime in which an infinite component exists, the so-called {\em percolation regime}, is a regime in which data can in principle be transmitted over long (infinite) distances via multiple hops. 

In the past 60 years since the first proof of continuum percolation in the Poisson--Gilbert graph, this field has expanded tremendously and presence as well as absence of percolation has been established in a great variety of static spatial models that extend the classical setting in various ways. Keeping the application area of D2D communication systems in mind, we highlight percolation results in so-called signal-to-interference-to-noise-ratio graphs~\cite{sinrPerc,tobias2020signal,jahnel2022sinr} (SINR graphs), where the existence of an edge not only depends on the mutual distance but also on the density of other devices in the vicinity. Notably, in SINR graphs one can observe an in-and-out of percolation in the intensity parameter, since for large intensities, the interference reduces the connectivity. On the other hand, let us mention Cox-percolation results~\cite{hirsch2019continuum,jahnel2022phase,hirsch2022sharp}, where the underlying point cloud is a Poisson point process with a random intensity measure that can for example be used to model urban environments such as street systems. Going further in this direction, in~\cite{HJCW22} sub- and supercritical regimes of percolation are established in a static spatial random graph model where the edge-drawing mechanism is designed to reflect realistic features of a connectivity network of moving devices in an urban environment. In short terms, any pair of devices, given by a Cox point process in which the environment is a planar random segment process, is connected by an edge iff, on their path through the environment, they spend enough time on the same segment in sufficiently close proximity. A variety of phase transitions in the different parameters can be observed, where again the speed parameter for the device movement, in certain settings, features an in-and-out of percolation. 

Going beyond the setting of static random graphs and their percolation behavior, the next natural step is to investigate data transmission on random graphs. Here the starting point is to consider growth processes in which a message is present at initial time at some (typical) vertex and then crosses the edges of the graph according to some passage times. This is the setting of {\em first-passage percolation} (FPP)~\cite{auffinger201750,kesten2003first} and results often concern the speed of the data transmission over large (Euclidean) distances and the associated (asymptotic) shape of the set of vertices that have received the message up to some fixed time. In~\cite{coletti2021limiting} such a shape theorem is presented for FPP on the supercritical Poisson--Gilbert graph. Conceiving the message as some infection, these growth models are often presented in the context of spatial probabilistic epidemiology and then it is natural to generalize the pure growth process of FPP to processes in which infected devices can also spontaneously heal (or messages can be dropped), giving rise to {\em contact processes} on random graphs~\cite{Li85,Li13}. In this setting the main concern usually becomes to show survival and extinction of the infection and results are available also for contact processes on random graphs in the continuum~\cite{menard2016percolation,bhamidi2021survival,nguyen2022subcritical}. In the context of D2D networks however, another antagonist of the pure growth mechanism is important. Regarding the infection as some malware, there can be special devices in the network that have the property to remove the malware from infected devices, but are otherwise indistinguishable from the infected or susceptible devices. This gives rise to {\em chase-escape processes}~\cite{durrett2020coexistence,tang2018phase,bernstein2022chase,beckman2021chase,hernandez2022distance} in which malware is transmitted from infected to susceptible devices and infected devices are patched. Again, main results are concerned with exhibiting regimes of survival and extinction of the malware in the system and such regimes have been established also in the continuum~\cite{hinsen2020phase,hinsen2020malware}. 

In the present manuscript we present two sets of results for chase-escape dynamics on supercritical percolation models in the continuum. In the first set of results we use the construction of~\cite{HJCW22}, as mentioned above, and consider a system of moving particles in an environment of line-segments in $\R^2$. However, we now go substantially beyond the setup in~\cite{HJCW22} and not only consider the resulting connectivity graph but introduce a chase-escape dynamics that respects the mobility and connectivity behavior of the network components. Controlling the long-range dependencies we prove results on survival and extinction of malware in certain parameter regimes and in particular feature again a phenomenon of in-and-out of survival with respect to the speed parameter. In the second set of results we complement these findings and present similar results in hybrid models in which the infection and patching processes are defined as an independent process on top of the (static) percolation model from~\cite{HJCW22}. 

The paper is organized as follows. In the remainder of this section we reproduce and extend our connectivity model from~\cite{HJCW22}. That is, the first network layer, described in Section \ref{route_sec}, is given by the street system. In Section~\ref{dev_sec} we describe the system of initial device positions on the street system and we introduce the paradigmatic mobility model for the devices given by the random waypoint model. In Section~\ref{con_sec} we introduce our notion of connectivity in the system and in Section~\ref{well_sec} establish local connectedness in the sense that almost-surely all devices have a finite degree. In Section~\ref{sec_results_WKP} we introduce the chase-escape dynamics and exhibit our main results for the global survival and absence of global survival of the malware including a theorem that establishes an in-and-out of survival with respect to the speed parameter. Next, Section~\ref{sec_results_FPP} features the hybrid model where the chase-escape dynamics is defined as an independent process on the connectivity graph and we present statements about global survival and absence of global survival under a variety of parameter regimes. Finally, in Section~\ref{sec_proofs} we present the proofs. 

\subsection{Street systems}\label{route_sec}
\begin{figure}[!htpb]
\centering
\input{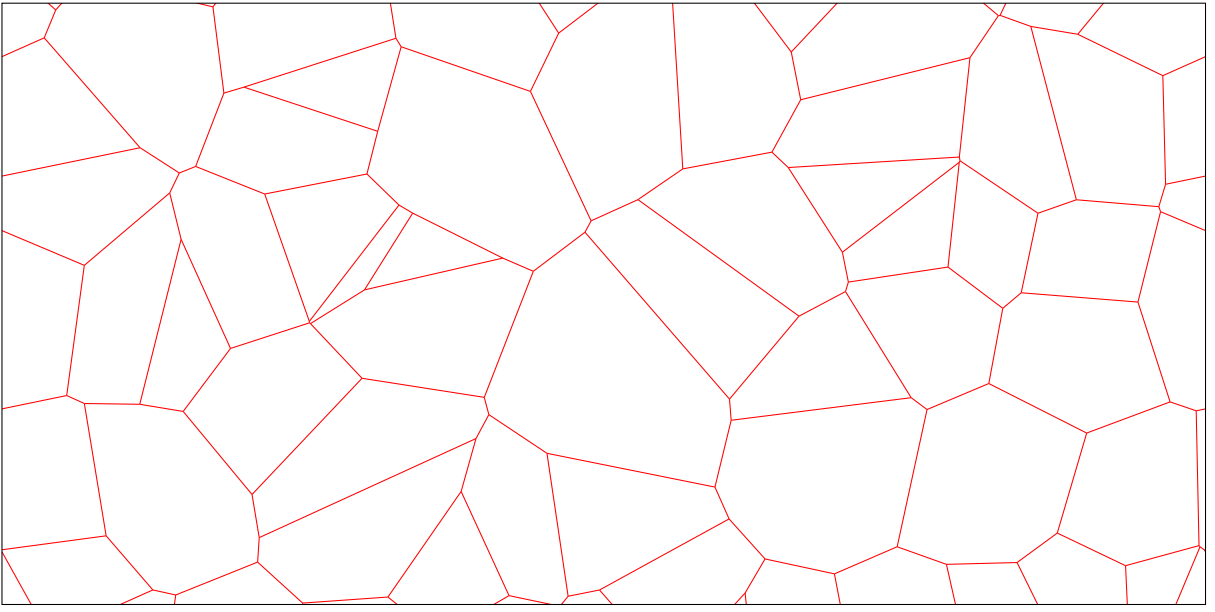}
\caption{Realization of a street system given by a Poisson--Voronoi tessellation.}
\label{Fig_Streets}
\end{figure}

As a first layer of the system, we consider a {\em stationary random planar segment process} 
$$S=\bigcup_{i\in I} \{S_i\}\subset \R^2$$
where each $S_i$ is a line segment in $\R^2$ of finite length. The line segments may intersect and we call each maximal unintersected interval in $S$ a {\em street}. The endpoints of each street are called {\em crossings} (even if they are dead-ends). In order to have something specific in mind, $S$ can could be for example some planar tessellation process like the {\em Poisson--Voronoi tessellation (PVT)} or the {\em Poisson--Delaunay tessellation (PDT)}, see Figure~\ref{Fig_Streets} for an illustration. 
We think of $S$ as to be an urban street system and note that for example the PVT indeed shares some common characteristics with medieval city topologies, see~\cite{courtat2012promenade}. 

In order to control spatial dependencies in $S$, we will almost exclusively work with systems that satisfy a quantitative mixing property called {\em stabilization}~\cite{stab1, stab2, stab3}. Additionally, we will assume that the system is connected in large regions with high probability and call this {\em asymptotic essential connectedness}~\cite{hirsch2019continuum,jahnel2022phase}. More precisely, let us 
define the square with side length $n\ge1$ centered at $x \in \R^2$ by
$$Q_n(x) = x + [-n/2, n/2]^2,$$ 
and abbreviate $Q_n = Q_n(o)$ and denote by $\dist(\varphi, \psi) = \inf\{|x-y|:\, x\in\varphi, y\in\psi\}$ the distance between sets $\varphi, \psi\subset\R^2$.
\begin{definition}[Stabilization and asymptotic essential connectedness]	\label{stabDef}
A stationary random segment process $S$ is called \emph{stabilizing} if there exists a random field of \textit{stabilization radii} $R = \{R_x\}_{x\in\R^2}$ that is measurable with respect to $S$ and $(S, R)$ is jointly stationary. Moreover, for $R(Q_n)=\sup_{y\in Q_n \cap\, \Q^2}R_y$ we have $\lim_{n\uparrow\infty}\P(R(Q_n) < n) = 1$ and for all $n \ge 1$, the random variables
 $$\Big\{f(S\cap Q_n(x))\one\{R(Q_n(x)) < n\}\Big\}_{x \in \varphi}$$
are independent for all non-negative bounded measurable functions $f$ and finite $\varphi \subset \R^2$ with $\dist(x, \varphi \setminus \{x\}) > 3n$ for all $x \in \varphi$.

Furthermore, $S$ is called {\em asymptotically essentially connected} if it is stabilizing and for all sufficiently large $n \ge 1$, whenever $R(Q_{2n})<n/2$, we have that $|S\cap Q_n|>0$ and $S\cap Q_n$ is contained in one of the connected components of $S\cap Q_{2n}$.
\end{definition}
We note that PVTs and PDTs are stabilizing~\cite{hirsch2019continuum} but Manhattan grids or (rectangular) Poisson-line tessellations are not stabilizing.

Finally, in what follows, we will always assume that the expected number of streets and the expected number of crossings in the unit square is finite and that the intensity of the street system is normalized, i.e., $\E[|S\cap Q|]=1$, where $|\cdot |$ denotes the total edge length.

\subsection{Mobile devices}\label{dev_sec}
\begin{figure}[!htpb]
\centering
\input{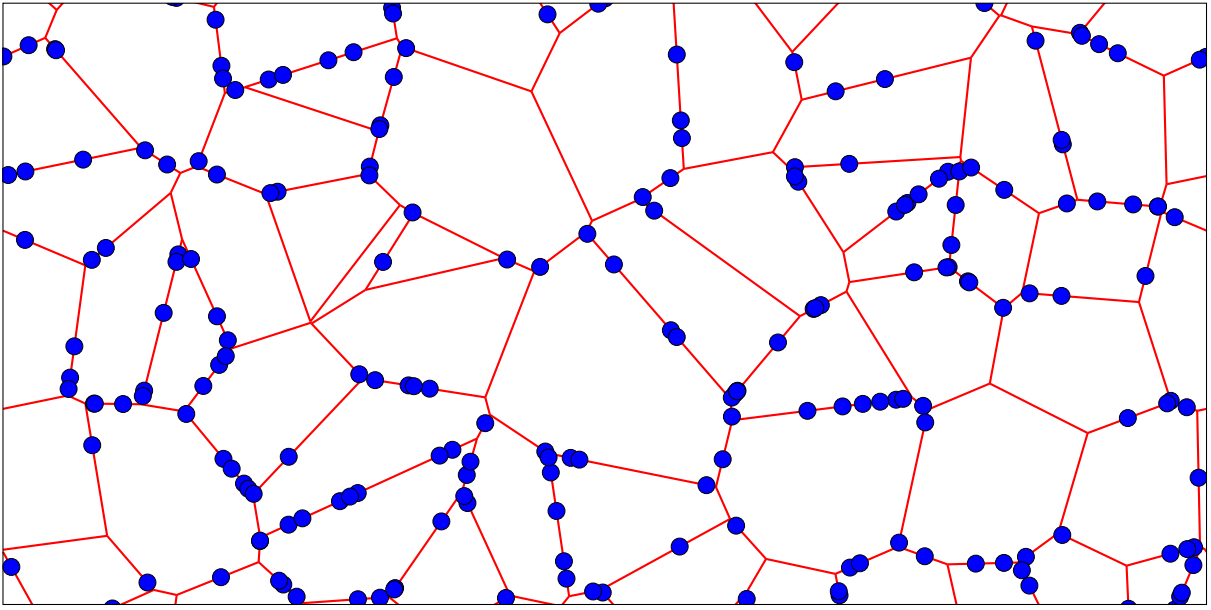}
\caption{Realization of initial device positions (blue) confined to a street system given by a Poisson--Voronoi tessellation.}
\label{Fig_Urban_Points}
\end{figure}

As a second layer of our model we represent {\em initial positions of  devices} via a {\em Cox point processes} $X^\la = \{X_i\}_{i\ge1}$. That is, $X^\la$ is a Poisson point process with a random intensity measure
$$\L_{S}(A) = \la |S\cap A|,\qquad\qquad A\subset \R^2 \text{ measurable},$$
where $\la\ge 0$ is a parameter that represents the expected number of devices per unit of street length, see Figure~\ref{Fig_Urban_Points}. 
Note that if $S$ is stationary so is $X^\la$.

Having defined the initial positions of the devices via a Cox point process on the street system, we now present our model for the device mobility. We start by considering the probability kernel
$$
\k^S(x, {\rm d} y),\qquad\qquad x\in S,
$$
which depend on $S$, and assume that the {\em support} of $\k^S(x,\d y)$, defined via
$$\supp(\k^S(x,\d y))=\{y\in S\colon \k^S(x,B_\e(y))>0\text{ for all }\e>0\},$$
is contained in $S$ for any $x\in S$. We say that $\k$ has (uniformly) {\em bounded support} if there exists $K>0$ such
that $$|\supp(\k^S(x,\d y))|<K\qquad\qquad \text{almost all }S,$$
where, with a slight abuse of notation $|\cdot|$ denotes the Lebesgue measure on $\R^2$. 
Moreover, we will always assume that $\k$ is {\em translation covariant}, i.e.,
$$\k^S(x,A)=\k^{S-x}(o,A-x)\qquad\qquad \text{almost all }S, \text{ measurable }A\subset \R^2,\,   x\in \R^2.$$
The kernel $\k^S(x,\d y)$ serves as a {\em waypoint kernel}~\cite{bettHart} and an example is given by the uniform distribution on $S\cap B_L(x)$, see~\cite{HJCW22}, where $B_L(x)\subset\R^2$ denotes the disc of radius $L$ centered at $x\in \R^2$.
Furthermore, we require that the kernel is {\em finitely dependent}, i.e., that there is a constant $H>0$ such that for almost-all $S$ and $x\in S$, $\k^S(x,\d y)=\k^{S'}(x,\d y)$ for all $S'$ with $S\cap B_H(x)=S'\cap B_H(x)$.

It will become useful to consider waypoint kernels that allow for very small displacements. In order to formalize this, we call the kernel {\em $c$-well behaved} if
$$
B_c(x) \cap S \subset \text{supp}\big(\k^S(x, \d y)\big)\qquad\qquad\text{almost all }S,\, x\in S,
$$
where $\k^S(x,y)$ is the density 
and say that it is {\em well behaved} if it is $c$-well behaved for some $c > 0$.

\begin{figure}[!htpb]
\centering
\input{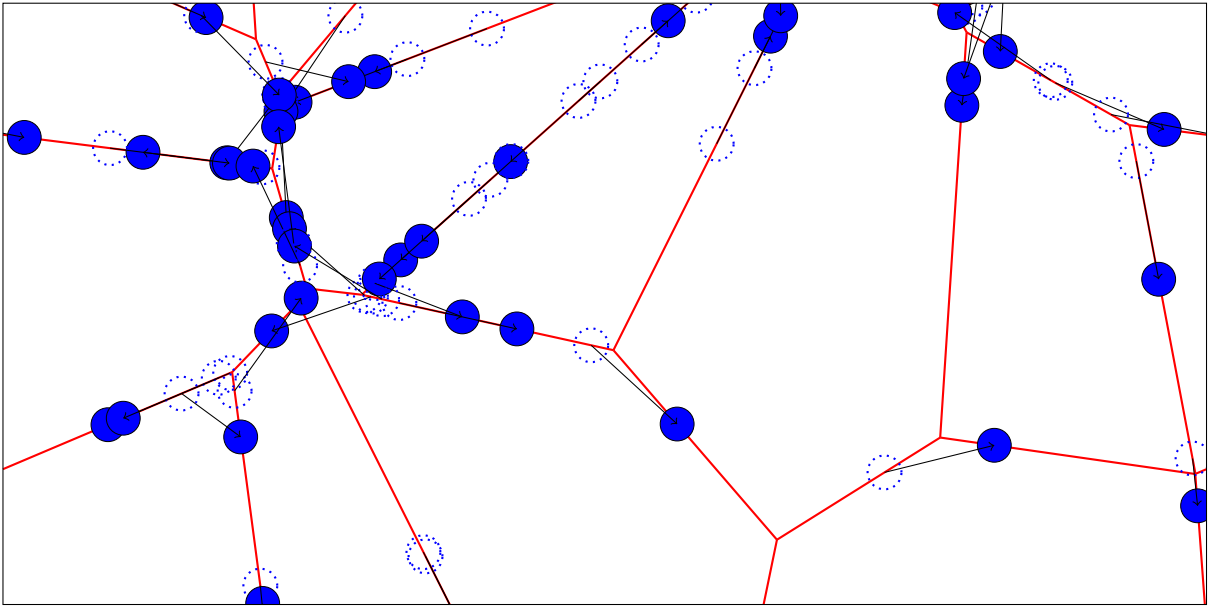}
\caption{Realization of initial device positions (dotted blue) confined to a street system given by a Poisson--Voronoi tessellation and their respective positions at a fixed positive time (blue), with arrows indicating the corresponding displacement.}
\label{Fig_Mobil_Points}
\end{figure}

Then, each device $X_i\in X^\la$ picks independently a {\em target location} $Y_i\in S$ according to $\k^S(X_i,\d y)$. Further, we equip each device $X_i\in X^\la$ with an i.i.d. velocity $V_i$, drawn from the distribution $\muv$ for which we assume that $0<v_{\rm min}\le v_{\rm max}<\infty$, where $v_{\rm min}:=\sup\{v\colon \muv([v,\infty))=1\}$ and $v_{\rm max}:=\inf\{v\colon \muv([0,v])=1\}$. Then, the path of the device $X_i$ is defined iteratively as follows. $X_i$ moves with speed $V_i$ to $Y_i$ along the shortest route in $S$ that connects $X_i$ and $Y_i$. If there are multiple shortest routes, the devices chooses one such route independently and uniformly at random. It then immediately returns to its starting position with the same velocity along the same shortest path. This procedure is then repeated indefinitely, see Figure~\ref{Fig_Mobil_Points} for an illustration.  

\subsection{Connectivity of mobile devices}\label{con_sec}
By a slight abuse of notation, we also denote by $X_i=(X_{i,t})_{t\ge 0}$ the trajectory of device $X_i$ in $S$. For any pair of devices $X_i$ and $X_j$ we then consider the set of {\em contact times}
$$
Z(X_i,X_j)=\{t\ge 0\colon |X_{i,t}-X_{j,t}|<r\text{ and } X_{i,t}, X_{j,t}\text{ are on the same street}\},
$$
as the times where the devices are on the same street and $r-$ close together, where $r>0$ is another parameter in the model, the {\em connectivity threshold}. 
Let us note that the constraint that contact times require the devices to be on the same street can be seen as a very strict {\em shadowing} assumption~\cite{BEEG21}. 

\begin{figure}[!htpb]
\centering
\hspace{-7pt}
\includegraphics[scale = 0.53]{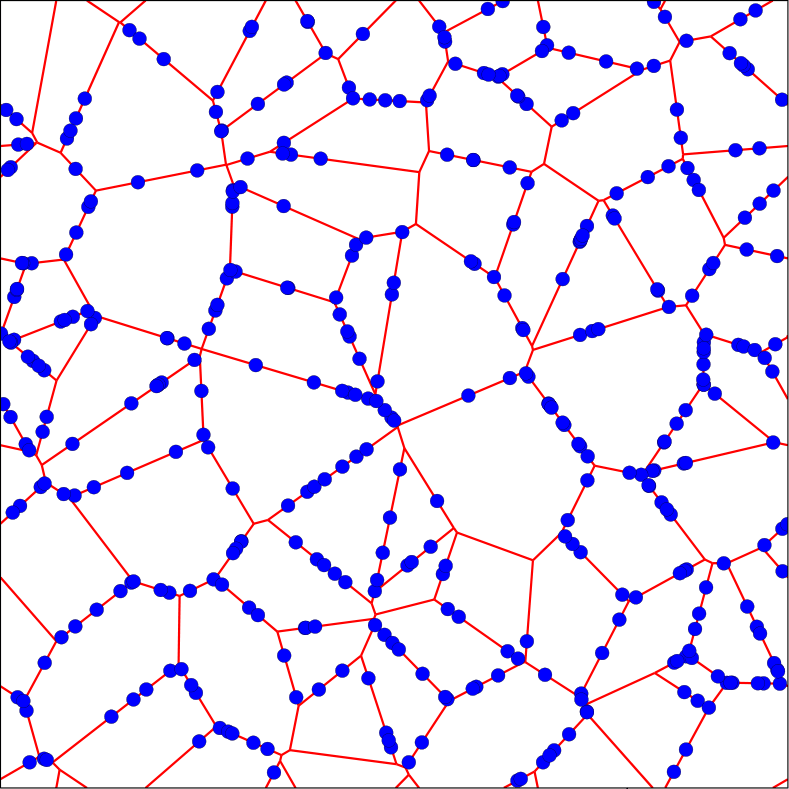}
\includegraphics[scale = 0.53]{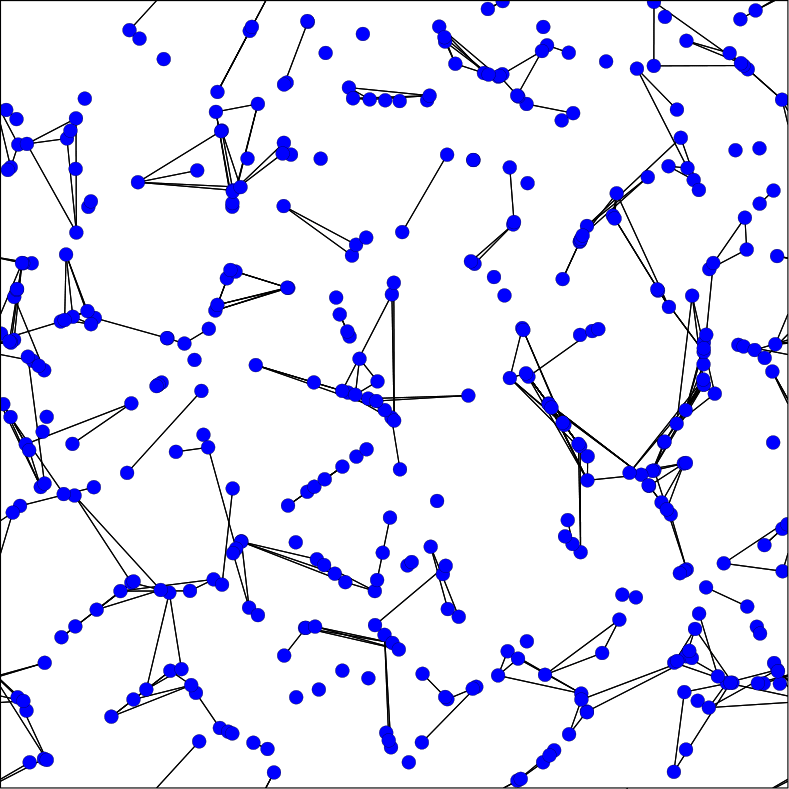}
\includegraphics[scale = 0.53]{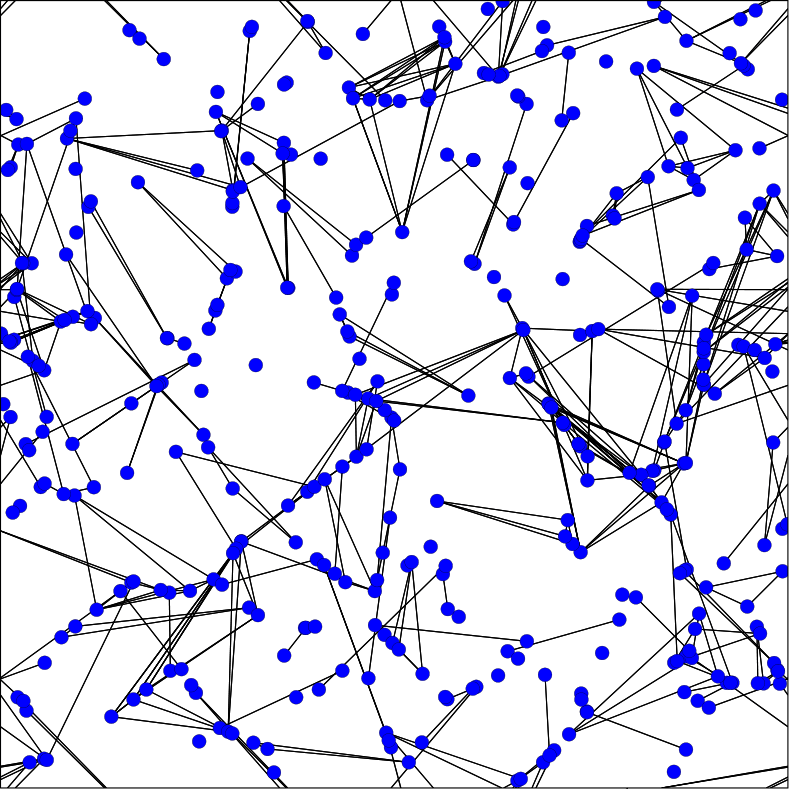}

\vspace{2pt}
\includegraphics[scale = 0.53]{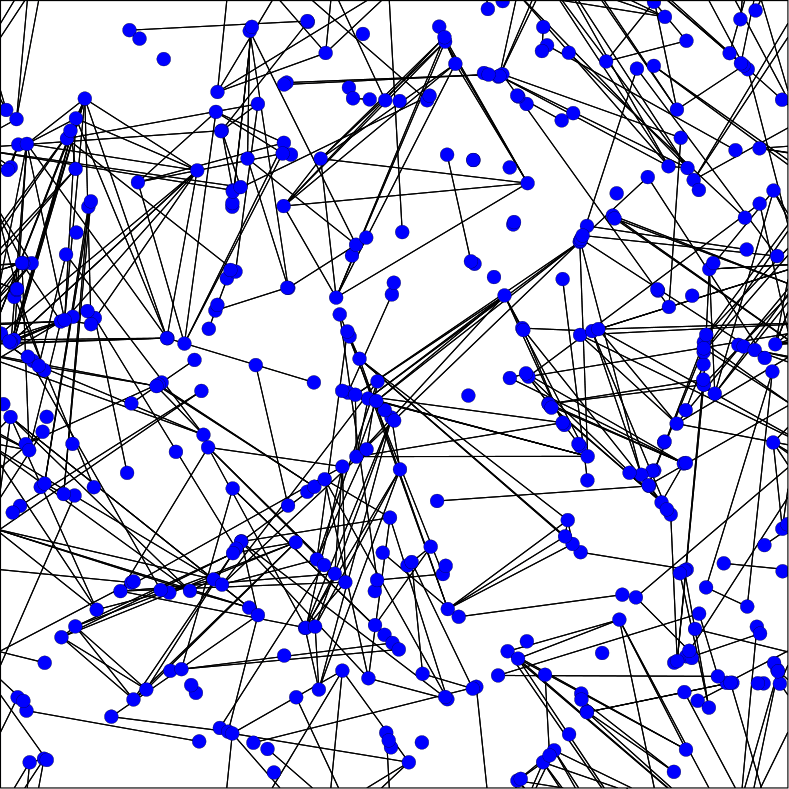}
\includegraphics[scale = 0.53]{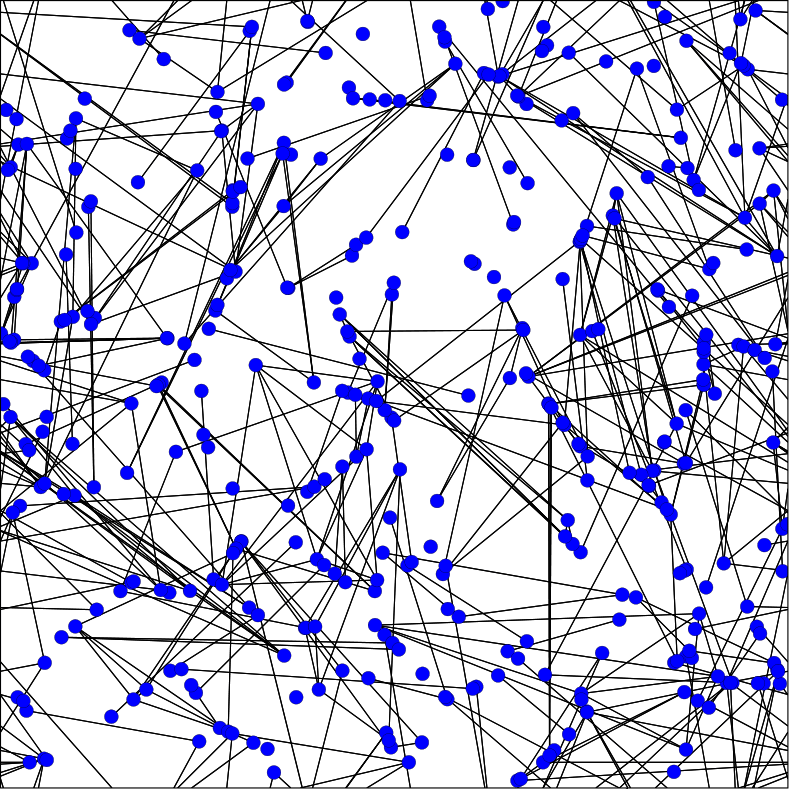}
\includegraphics[scale = 0.53]{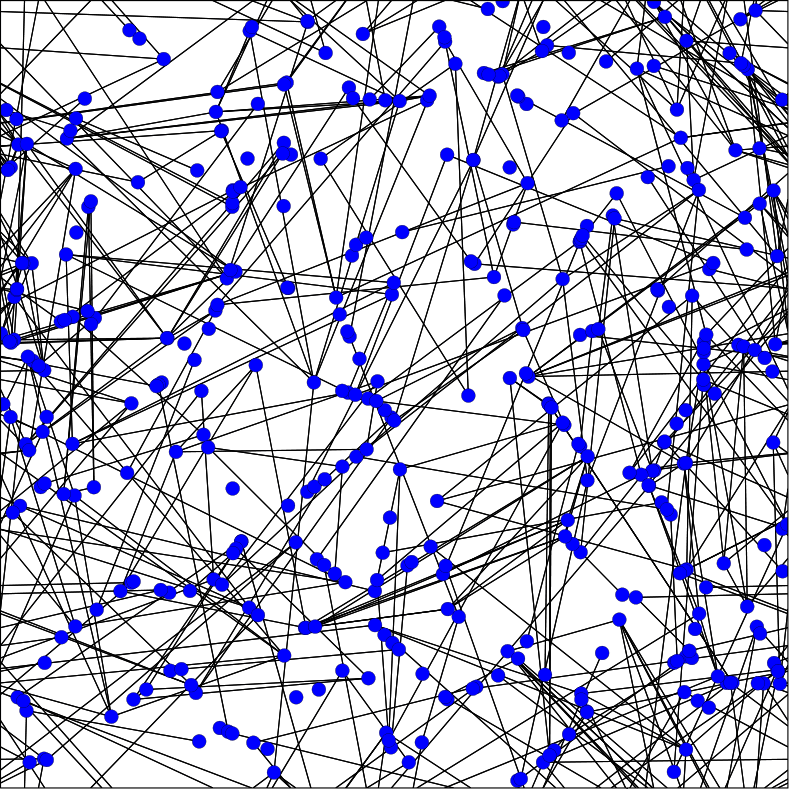}
\caption{Realization of initial device positions (blue) on a street system given by a Poisson--Voronoi tessellation (red, only drawn in the initial picture, visually suppressed otherwise). Between two devices $X_i$ and $X_j$ an edge is drawn if $t\in Z(X_i,X_j)$ for $t=0,1,\dots,5$. We highlight that the total number of edges in this sequence of realizations stays roughly constant, however the edges become longer. This is due to the fact that we plot initial positions, but devices move and hence, over time, connect to more distant devices.
}
\label{Fig5}
\end{figure}

Next, we assume that any pair of devices $X_i$ and $X_j$ has an iid {\em infection time} $\rho(X_i,X_j)$, drawn from the distribution $\varrho$ on $(0,\infty)$. Then, we say that a {\em transmission from $X_i$ to $X_j$ is possible at time $t\ge 0$} if
$$[t-\rho(X_i,X_j),t]\subset Z(X_i,X_j),$$
see Figure~\ref{Fig5} for an illustration of the resulting space-time connectivity network.

Let us note that the introduction of random transition times is a substantial generalization compared to the setup in~\cite{HJCW22}, where only the case $\varrho=\delta_\rho$ is considered for some fixed $\rho\ge0$. 

\subsection{Degree of devices}\label{well_sec}
In order to avoid scenarios in which a single device can transmit to infinitely many devices in the infinite time horizon that we consider, let us finally introduce a convenient condition for the {\em local connectedness} of the connectivity graph. For this, let $\ell_S(x,y)\subset \R^2$ be the shortest path between $x$ and $y$ on $S$ and write 
$$\ell_S(x)=\sup\{|\ell_S(x,y)|\colon y\in\supp(\k^S(x,\d y))\}+r/2$$
for the length of the shortest path starting in $x$ towards any reachable target on $S$ plus half the connectivity range. Then, consider the following {\em Cox--Boolean model with geostatistical markings} in which $X_i,X_j\in X^\la$ are connected whenever
$$
|X_i-X_j|\le \ell_S(X_i)+\ell_S(X_j)
$$
and note that this can be seen as a Boolean model in which the (random) discs, associated to each Cox point, depend on the underlying environment $S$. 
In particular, if $X_i$ can transmit to $X_j$, there also exists an (undirected) edge between $X_i$ and $X_j$ in the Cox--Boolean model with geostatistical markings. Using this, for $X_i\in X^\la$, we define its degree 
$$\deg(X_i)= \#\{X_j \in X^\la \setminus X_i \colon |X_i-X_j|\le \ell_S(X_i)+\ell_S(X_j)\},$$
which serves as an upper bound for the degree of $X_i$ in the original model.  
Now, we call the model {\em locally connected} if 
\begin{align}\label{locconn}
\P(\exists X_i\in X^\la\text{ such that }\deg (X_i)=\infty)=0
\end{align}
and, in the sequel, we will always assume that our system is locally connected. Let us also mention that, in order to ease notation, we use generic symbols $\P$ and $\E$ for the distribution and expectation of our model, even under changing parameters. 

The following proposition establishes a condition under which the network is locally connected. 
\begin{proposition}[Local connectedness]\label{locconn_def}
If $\k$ has bounded support, $S$ is exponentially stabilizing and $|S\cap Q_1|$ has exponential moments, then the network is locally connected.
\end{proposition}
The proof rests on a contradiction that follows from a first-moment method for the expected degree of a typical point using the good control on stabilization and moment properties of the street system. We present the details in Section~\ref{proof_locconn}. Let us note that the conditions are satisfied for our standard examples for street systems given by PVT and PDT, since they are exponentially stabilizing~\cite{hirsch2019continuum} and have exponential moments, see~\cite{JaTo19}. Also, as can be seen from the proof, the moment condition can be substantially relaxed.  

Having defined now our system of moving devices in an urban street system and their connectivity relations, in the following section, we introduce a {\em malware} into the network that propagates like an {\em infection} starting from a single typical device.   

\section{Chase-escape on the dynamic graph}\label{sec_results_WKP}
In the final model layer, we introduce an infection into the system, which is carried by one device $X_o$ that is chosen at random and will be referred to in the sequel as the {\em typical device}. The remaining devices are, at time zero, {\em susceptible} to the infection. As a counter measure we introduce an independent (Cox) point process of so-called {\em white knights}, which is, at time zero, a homogeneous Poisson point process $Y^{\laW}$ with linear intensity $\laW\ge 0$ conditioned on $S$ and following the same mobility scheme as devices in $X^\la$. From a modeling perspective, we are here guided by the idea that white knights are also regular devices and their general behavior is indistinguishable from other devices. However, white knights are not susceptible to the infection. Even more, each white knight carries a {\em patching} software that allows it to remove the infection from any device that wants to infect it with malware and update the defensive mechanism on the attacking device. Let us highlight, that there is no interaction between susceptible devices and white knights. Moreover, the process of susceptible devices is decreasing and the process of white knights is increasing, while the process of infected devices is non-monotone. More formally, let us introduce the {\em Palm version of susceptible devices} via
\begin{align}\label{Palm}
    \E^*[f(X^\la,Y^{\laW})]=\frac{1}{\la}\E\Big[\sum_{X_i\in X^\la\cap Q_1 }f(X^\la-X_i,Y^{\laW}-X_i)\Big]\qquad\text{ bounded, measurable testfunctions } f,
\end{align}
and note that under $\P^*$, there exists a device at the origin with probability one. We call this the typical device $X_o$ and assume that it is {\em infected} at time zero. 

\begin{figure}[!htpb]
\centering
\hspace{-7pt}
\includegraphics[scale = 1.15]{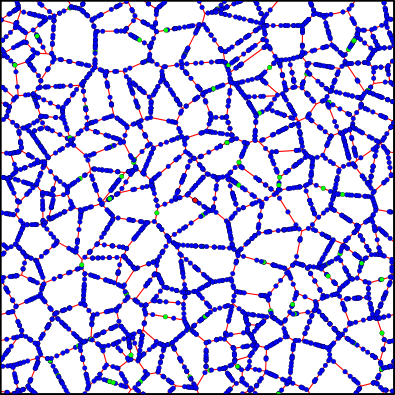}
\includegraphics[scale = 1.15]{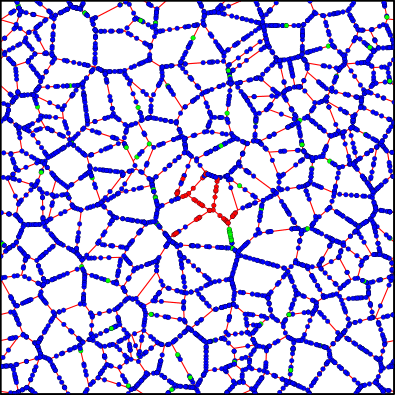}
\includegraphics[scale = 1.15]{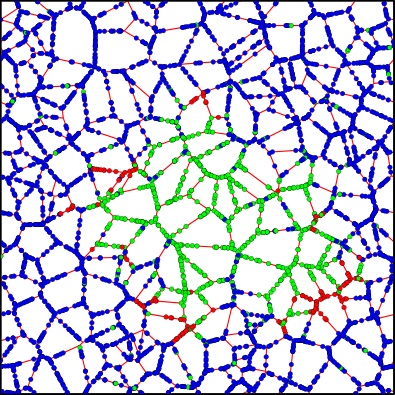}
\caption{Propagation of infected devices (red) on a street system given by a Poisson--Voronoi tessellation. In the initial state (left) there is exactly one infected device present in the center and the remaining devices are either susceptible (blue) or white knights (green). At some small positive time (middle) further devices in the vicinity of the initially infected device have become infected and have started to make contact to white knights. At some later time (right) infected devices are only present along the boundary of the set of white knights in the center region. 
}\label{Fig6}
\end{figure}

In order to describe the infection and patching mechanism, instead of a single infection time distribution $\varrho$, we consider two distributions $\varrho_{\rm I}$, respectively $\varrho_{\rm W}$, for the {\em infection times}, respectively the {\em patching times}, and assume that they both have a density with respect to the Lebesgue measure on $(0,\infty)$. We denote by $\varrhoImin= \inf\{x\ge 0\colon x \in \text{supp}(\varrho_{\rm I})\}$ the minimal infection time and analogously by $\varrhoWmin$ the  minimal patching time. Now, the infection will be transmitted from an infected device $X_i$ to a susceptible device $X_j$ after completion of the first infection time $\rhoI(X_i,X_j)$, drawn from $\varrho_{\rm I}$. 
Analogously, the patch will be transmitted from a white-knight device $X_i$ to an infected device $X_j$ after completion of the first infection time $\rhoW(X_i,X_j)$, drawn from $\varrho_{\rm W}$. In order to formalize the whole process, let $S_t$, $I_t$, and $W_t$ denote the sets of susceptible, infected and white-knight devices at time $t\ge 0$, where $S_0=X^\la\setminus X_o$, $I_0=X_o$, and $W_0=Y^{\laW}$. In particular, 
$\big((I_s)_{0\le s<t},(W_s)_{0\le s<t}\big)$ fully describes the state of the system up to time $t$. Then, at time $t$, we add new infected devices 
\begin{align*}
    I_t\setminus \bigcup_{0\le s<t} I_s=\Big\{X_i\in \bigcup_{0\le s<t} S_s\colon \exists X_j \in \bigcap_{s \in [t-\rhoI(X_j,X_i),t)}I_s  \text{ such that }   [t-\rhoI(X_j,X_i),t) \subset Z(X_j,X_i) \Big\},
\end{align*}
and then add also new white knights
\begin{align*}
    W_t\setminus \bigcup_{0\le s<t} W_s=\Big\{X_i\in X^\la\colon& \exists X_j \in W_{t-\rhoW(X_j,X_i)} \text{ such that }  \\
    &[t-\rho_{\rm W}(X_j,X_i),t) \subset Z(X_j,X_i)\text{ and } X_i \in I_{t-\rho_{\rm W}(X_i,X_j)}\setminus \bigcup_{0\le s<t} W_s\Big\}.
\end{align*}
With this construction the process is well-defined since $\varrho_{\rm I}$ and $\varrho_{\rm W}$ put no mass on zero and the number of devices that is infected is finite for any finite time by the local connectedness condition. Note that indirectly we have implemented a tie-breaker rule that favors the infection. This can be seen by considering a situation where a device $X_j$ becomes a white knight at time $t$, however it can still finish a transmission of the infection at time $t$. We present a realization of the chase-escape dynamics in Figure~\ref{Fig6}.

We say that the infection {\em survives locally} if $|I_t|>0$ for all $t\ge 0$, it {\em survives globally} if $|I_t|>0$ for all $t\ge 0$ and $|\bigcup_{t \ge 0} I_t|=\infty$, and the infection {\em goes extinct} if $|I_t|=0$  for some $t\ge 0$.
Note that, even if the street system is connected, with positive probability, there can be isolated sets of devices that are not able to form connections to the rest of the devices. If the typical device $X_o$ is contained in such a cluster and the cluster contains no white knight, the infection survives indefinitely, but no global outbreak can be observed.

\medskip
Our goal is to isolate regimes of global survival and extinction of the infection using the transmission mechanism we just described. 
Let us start with regimes that ensure global survival of the infection even in the presence of white knights. For this, we guarantee the existence of a sequence of streets on which the infection can be transmitted to infinity with positive probability. 

Note that, due to the transmission mechanism, streets need a certain minimal length in order to be useful for transmissions. 
We thus define the thinned system $S^a=\{s\in S \colon |s|\ge a\}$ that only contains streets of length greater or equal than $a\ge 0$. Let us further denote by $R^a_x$ the distance of the furthest point from $x\in \R^2$ that is reachable without crossing $S^a$, i.e., we define the cell of $x$ by
$$C^a_x=\sup\{y \in \R^2\colon \text{ there exists a continuous function $f\colon [0,1] \mapsto \R^2 \setminus S^a$ with }f(0)=x, f(1)=y\}$$
and define $R^a_x=\sup\{|y-x|\colon y \in C^a_x\}$. We say that the thinned graph is {\em $R^a$-connected} if $\lim_{n\uparrow \infty}\P\big(\sup_{x \in Q_n\cap \Q^2}R^a_x<n\big) =1$ and define 
\begin{align}\label{def_ac}
a_{\rm c}:=\sup\{a>0\colon S^a \text{ is } R^a \text{-connected}\}.
\end{align}
In words, for $a<a_{\rm c}$, the thinned street system satisfies a slightly stronger version of asymptotically essentially connectedness, and we have some control over the cell sizes. Note that  
$a_{\rm c}>0$ for asymptotically connected street systems, see~\cite[Theorem 2.8]{HJCW22}. 

In order to exhibit a regime of global survival, we will furthermore require a slightly stronger version of well behavedness. We say that a kernel $\k$ is {\em $c$-continuous} if, for almost-all $S$ and $x\in S$, the uniquely defined absolutely-continuous part of $\k^S(x,\d y)$ is non trivial and its density $\k^S(x,y)$ with respect to the Lebesgue measure on $S$ satisfies
\begin{align*}
    \inf_{y \in B_c(x)}\k^S(x,y)>0.
\end{align*}
Note that $c$-continuous implies $c$-well behaved.

Using these definitions, we can state our first main result on global survival of the infection. 

\begin{theorem}[Global survival]\label{Thm_Survival}
Let $\k$ be $c$-continuous for some $c>0$ and $\varrhoWmin>\varrhoImin>0$. Then, for all sufficiently small $\vmin$ such that $0<\vmin \varrhoImin < \min(a_{\rm c}/2,r,c/2) $, there exists $\lac$ such that for all $\la>\lac$ and all $\laW\ge 0$,
$$\P(\text{infection survives globally})>0.$$
\end{theorem}
The punchline of the preceding statement is the fact that the critical device intensity is independent of the white-knight intensity. In other words, no matter how dense we can pack white knights into the system, survival is guaranteed as long as we have sufficiently many devices in the system. Roughly, the reason for this phenomenon is that the white knights simply act too slow since $\varrhoWmin>\varrhoImin>0$ when there are sufficiently many devices that can transmit the infection to their neighbors before white knights are able to eliminate it. The idea of the proof is to establish an infinite cluster of good streets that have the property that, if the malware enters the street on one side, it must also exit the street on the other side. We present the details in Section~\ref{ProofThm_Survival}.

\medskip
Next, we present our result on the extinction of the infection in infinite components.
\begin{theorem}[Extinction]\label{Thm_Extinction}
Consider a street system given by a Poisson--Voronoi or Poisson--Delaunay tessellation. Let $\k$ be well behaved with bounded support and assume that $r>\vmax\varrhoWmin$ and $\varrhoImin>\varrhoWmin>0$. Then, for all $\lambda\ge 0$, there exists $\laWc$ such that for all $\laW>\laWc$ we have that
$$\P(\text{infection survives globally})=0.$$
\end{theorem}
The proof rests on a multi-scale argument and will be given in Section~\ref{Proof_Thm_Extinction}.

\medskip
Finally, we present a result that features the remarkable non-monotone behavior of our model with respect to the speed parameter. More precisely, the following result states that both, too high and too low velocities, are detrimental for the propagation of the infection. We state the result only for fixed velocities and infection and patch times and note that a corresponding result for general $\mu_{\rm v}$, $\varrho_{\rm I}$ and $\varrho_{\rm W}$ also holds.
\begin{theorem}[Speed-dependent survival and extinction]\label{thm_InOut}
Consider a street system given by a Poisson--Voronoi or Poisson--Delaunay tessellation. Let $\k$ have bounded support and be $c$-well behaved for some $c>0$. Assume further that $\varrho_{\rm W}=\delta_{\rho_{\rm W}}$ and $\varrho_{\rm I}=\delta_{\rhoI}$ with $\rho_{\rm W}>\rhoI$, and let $0<v_o< \min(a_{\rm c}/2,r,c/2)/\rhoI $. Then, there exists $\la_{\rm c}>0$ and $\laWc>0$ (independent of $\la_{\rm c}$ and $v_o$), such that for all $\la >\la_{\rm c}$ and all $\la_{\rm W}>\laWc$ we have that
\begin{enumerate}
\item[(1)] there exists $v_o>v_{\rm c}(\la,\la_{\rm W})>0$ such that for all $\muv=\delta_{v}$ with $v<v_{\rm c}(\la,\la_{\rm W})$ we have 
$$\P(\text{infection survives globally})=0,$$
\item[(2)] for $\muv=\delta_{v_o}$ we have $\P(\text{infection survives globally})>0$, and
\item[(3)] there exists $\infty>v^{\rm c}>v_o$ (independent of $\la$ and $\la_{\rm W}$), such that for all $v>v^{\rm c}$ we have
$$\P(\text{infection survives globally})=0.$$
\end{enumerate}
\end{theorem}
In words, the choice $\la >\la_{\rm c}$ puts us in a parameter regime in which survival of the infection is possible, as described in Theorem~\ref{Thm_Survival}. However, in the presence of sufficiently many white knights and sufficiently large or small velocity we are again in regimes of extinction. One important aspect here is that sufficiently small velocities and large white-knight intensities lead to absence of global survival also in the case where the infection is faster than the patching. We present the proof in Section~\ref{Proof_thm_InOut}.

\section{Chase-escape on the connectivity graph}\label{sec_results_FPP}
In the previous Section~\ref{sec_results_WKP} we exhibited rigorous results for chase-escape on the dynamic graph. One important feature here is that the device mobility is intimately linked to the transmission mechanism, both of the infection as well as the patching, leading to a highly correlated non-Markovian space-time process. As a consequence, currently it seems out of reach to derive approximate values for the various critical parameters predicted by our rigorous analysis, in ways other than computer simulations as presented for example in~\cite{MultiAgent,MultiAgent_AAMAS}. As it turns out, for this a major part of the computational effort has to go into simulating the movement of the devices in the environment. In an effort to produce large numbers of sample paths of the chase-escape dynamics it is therefore desirable to uncouple the connectivity graph and the epidemiological process and derive corresponding results in this {\em mean-field type} setting. This is our primary goal in this section. 

To start, recall the connectivity rule, introduced in Section~\ref{con_sec}, where for every pair of devices $X_i$, $X_j$ the set of connection times is denoted $Z(X_i,X_j)$. Now, we want to fix the infection time $\rho \ge 0$, i.e., assume that $\varrho=\delta_\rho$ and define the \emph{connectivity graph} 
$$g_{\rho,r}(S,X^\la \cup Y^{\laW}),$$
in which any two vertices $X_i, X_j \in X^\la\cup Y^{\laW}$ are connected by an edge if there exists a $t\ge 0$ such that $[t,t+\rho]\subset Z(X_i,X_j)$.
Let us highlight that this connectivity graph is {\em static} even though the edge-drawing mechanism is based on an underlying space-time process. 
The main computational advantage of this connectivity graph lies in the fact that the connections can be determined purely from the device trajectories and speeds without the need to simulate in which connection interval the connection is realized. This avoids the computationally expensive collective movement of the devices.

Now, we consider the chase-escape model as an independent space-time process on realizations of the connectivity graph. This also brings us much closer to the existing literature on chase-escape dynamics on fixed and random graphs, see for example~\cite{durrett2020coexistence,tang2018phase,bernstein2022chase,beckman2021chase,hernandez2022distance,hinsen2020phase,hinsen2020malware}. More precisely, we consider \linebreak 
$g_{\rho,r}(S,X^\la \cup Y^{\laW})$
under the Palm distribution with respect to the $X^\la$ as described in~\eqref{Palm} and let $N(X_i)$ denote the set of direct neighbors of vertex $X_i$. 
Further, we consider two transmission-time distributions $\varrhoI$ and $\varrhoW$ on $[0,\infty)$ that govern the iid time $\rhoI(X_i,X_j)$, respectively $\rhoW(X_i,X_j)$, that is needed to pass an infection, respectively a patch, from device $X_i$ to $X_j$. Then, at time zero, the set of infected devices $I_0$ is given by the typical device $X_o$, which is positioned at the origin $o\in \R^2$ with probability one, i.e., $I_0=\{X_o\}$. For the initial set of white knights we have $W_0=Y^\laW$. Then, as in Section~\ref{sec_results_WKP}, we define the process iteratively by considering the system $\big((I_s)_{0\le s<t},(W_s)_{0\le s<t}\big)$ up to time $t>0$ and then define the newly infected devices at time $t$ by
\begin{align*}
    I_t\setminus \bigcup_{0\le s<t} I_s=\{X_i\in \bigcup_{0\le s<t} S_s\colon N(X_i) \cap \bigcap_{s \in [t-\rhoW(X_j,X_i),t)}I_s \neq \emptyset \},
\end{align*}
where $S_s=X^\la\setminus I_s$, and the newly patched devices at time $t$ by
\begin{align*}
    W_t\setminus \bigcup_{0\le s<t} W_s=\{X_i\in X^\la\colon \exists X_j \in W_{t-\rhoW(X_i,X_j)}\cap N(X_i), X_i \in I_{t-\rhoW(X_i,X_j)}\setminus \bigcup_{0\le s<t} W_s\}.
\end{align*}
Let us comment on the construction. As mentioned above, the model can be seen as an approximation of chase-escape on the dynamical graph where we associate to every edge an independent random time that represents the time that it takes to transmit the infection. The distribution of this transfer time can be, for example, sampled by considering the transmissions of a typical device towards its neighbors in an independent simulation step. In other words, we insert information about the typical transmission behavior and apply it independently to every edge of the static graph. Let us mention that this approach is not unlike the approach used in~\cite[Section V]{MultiAgent}, where the spatial graph was replaced by a sequence of edges which had the length distribution of a typical edge in the PVT. 

We are again interested in conditions under which global survival is possible or impossible. Let us start with the preliminary observation in the case where the transmission times are fixed, i.e., $\varrhoI=\delta_{T_{\rm I}}$ and $\varrhoW=\delta_{T_{\rm W}}$ for some $T_{\rm W},T_{\rm I}>0$. In this case, survival and extinction essentially depends on the ordering of $T_{\rm W}$ and $T_{\rm I}$. We say that the graph $g_{\rho,r}(S,X^\la)$ is in the {\em supercritical percolation regime} if there is an unbounded connected component and note that~\cite{HJCW22} exhibits non-trivial criteria that guarantee the existence of a supercritical percolation regime. 
\begin{proposition}\label{Prop_FPP_fixed_T}
Assume that $g_{\rho,r}(S,X^\la)$ is in the supercritical percolation regime. Then, 
\begin{align*}\P(\text{infection survives globally})\begin{cases} 
>0& \text{ for } T_{\rm I}\le T_{\rm W}, \laW\ge0\\
=0& \text{ for } T_{\rm I}>T_{\rm W},\laW>0
\end{cases}
\end{align*}
\end{proposition}
The proof is based on a robust descending-chain argument that in fact generalizes to arbitrary supercritical graphs and is presented in Section~\ref{sec_Prop4-1}. 

Let us now present our main results for general $\varrhoI$ and $\varrhoW$. In order to prove global survival of the infection we need to require stronger notions of percolation of the connectivity graph. We present two alternative conditions in Theorem~\ref{Thm_Exp_Sur}, respectively Theorem~\ref{Thm_Exp_Sur_2}. 
First, in Theorem~\ref{Thm_Exp_Sur} we require a parameter regime for the street system and waypoint kernel that guarantees a strongly percolating structure similar to Theorem \ref{Thm_Survival}. 
Recall the definition of $a_{\rm c}$ from~\eqref{def_ac} and let us denote for any $b>0$ by $\varrhoIb$ the shifted version of $\varrhoI$ defined via $\int f(x/b)\varrhoIb(\d x)=\int f(x)  \varrhoI(\d x)$ for all measurable functions with bounded support $f$.
\begin{theorem}[Mean-field survival I]\label{Thm_Exp_Sur} 
Let $\k$ be $c$-well behaved for some $c>0$, assume that $\vmin$ satisfies $0<\vmin \rho < \min(a_{\rm c}/2,r,c/2)$ and that $\varrhoWmin>0$. Then, there exists $\la_{\rm c}<\infty$ such that for all $\la >\la_{\rm c}$ and $\laW\ge 0$ there exists $b_{\rm c}<\infty$ such that for $\varrhoIb$ with $b>b_{\rm c}$ we have that 
$$\P(\text{infection survives globally})>0.$$
\end{theorem}
The proof rests on percolation arguments similar to the ones used for the proof of Theorem~\ref{Thm_Survival}. The details are presented in Section~\ref{proofThm_Exp_Sur}.

Next, we present an alternative result for global survival which requires that a subgraph of the connectivity graph is in the supercritical percolation regime. More precisely, we define the random subgraph $g^{M,p}_{\rho,r}(S,X^\la)\subset g_{\rho,r}(S,X^\la)$ via the following thinning procedure.
Let $(\textrm{Ber}_p(X_i))_{X_i \in X^\la}$ be an iid field of Bernoulli random variables with parameter $p$ and define the vertex set of $g^{M,p}_{\rho,r}(S,X^\la)$ by
\begin{align*}
V^{M,p} = \{X_i \in X^\la \colon  N(X_i)\le M \textrm{ and } \textrm{Ber}_p(X_i) = 1\},
\end{align*}
where $N(X_i)$ denotes the number of neighbors in the graph $g_{\rho,r}(S,X^\la, Y^{\laW})$ that includes also the white knights. In words, vertices with more than $M$ neighbors are eliminated and additionally each of the remaining vertices is kept independently with probability $p$. The edge set of $g^{M,p}_{\rho,r}(S,X^\la)$ is then inherited from the edge set of $g_{\rho,r}(S,X^\la)$, but only those that connect $X_i,X_j\in V^{M,p}$.
It is not hard to see that 
$$\lim_{M \uparrow \infty, p \uparrow 1}g^{M,p}_{\rho,r}(S,X^\la) = g_{\rho,r}(S,X^\la)$$
weakly, however for global observables such as percolation this is non trivial.
We have the following result.
\begin{theorem}[Mean-field survival II]\label{Thm_Exp_Sur_2}
Assume that there exist $M>0$ and $p>0$ such that $g^{M,p}_{\rho,r}(S,X^\la)$ is in the supercritical percolation regime. Further assume that $\varrhoI$ has no atom in $0$. Then, for all $\laW\ge 0$, there exists $b_{\rm c}<\infty$, such that for all $\varrhoIb$ with $b>b_{\rm c}$ we have that 
$$\P(\text{infection survives globally})>0.$$
\end{theorem}
The proof rests on the idea that, for sufficiently fast infection rate, already the process of nodes with a maximal degree, which also transmit the infection faster than being patched, is in the supercritical percolation regime. 
We present the details in Section~\ref{proofThm_Exp_Sur_2}.

We conclude this section with our result on the absence of global survival. 
\begin{theorem}[Mean-field extinction]\label{Thm_Exp_Ext}
Consider a street system given by a Poisson--Voronoi or Poisson--Delaunay tessellation and assume that $\varrhoImin>\varrhoWmin$.
Then, there exists $\laWc<\infty$, such that for all $\laW>\laWc$ we have that 
$$\P(\text{infection survives globally})=0.$$
\end{theorem}
The proof is based on the idea that, for sufficiently large intensity of white knights, any connection between susceptible devices is accompanied by many white knight connections. In this case, it becomes highly unlikely that the infection is passed on before being eliminated by the white knights. In order to make this precise, we have to leverage local independence properties of the propagation model. We present the details in Section~\ref{ProofThm_Exp_Ext}.


\section{Proofs}\label{sec_proofs}
\subsection{Proof Proposition~\ref{locconn_def}}\label{proof_locconn}
\begin{proof}[Proof Proposition~\ref{locconn_def}]
First, let us assume that $\P(\exists X_i\in X^\la\text{ such that }\deg (X_i)=\infty)>0$. Then, by continuity of measures, there exists $R>0$ such that $\P\big(\exists X_i\in X^\la\cap Q_{R}\text{ such that}\deg (X_i)=\infty\big)>0$ and hence, by translation invariance, also $\P\big(\exists X_i\in X^\la\cap Q_1\text{ such that}\deg (X_i)=\infty\big)>0$. But the last inequality implies that $\E[\sum_{X_i \in X^\la \cap Q_1}\deg(X_i)]= \infty$. However, in order to derive a contradiction, we can use the Mecke–-Slivnyak Theorem twice to see that
\begin{align*}
    \E\Big[\sum_{X_i \in X^\la \cap Q_1}\deg(X_i)\Big] 
    &\le \E\Big[\sum_{X_i \in X^\la \cap Q_1} \sum_{X_j \in X^\la \setminus X_i}\one\{|X_i-X_j| \le \ell_S(X_i) + \ell_S(X_j)\}\Big]
    \\
    &=\la^2\E\Big[\int_{S\cap Q_1} \d x\int_{S }\d y \one\{|x-y| \le \ell_S(x) + \ell_S(y)\}\Big].
\end{align*}
Next, we bound the maximal reach $\ell_S$ by the street length in a stabilized region. For this, let 
$$K:= \sup\{|x-y|\colon x\in S, y\in\supp(\k^S(x,\d y))\}<\infty$$ denote the maximal reach of the kernel $\k$, which exists due to the assumption of a uniformly bounded support of $\k$. Using this, we can bound the length of any shortest path starting in $x\in S$ to any $y \in S\cap B_K(x)$, by the total street length of the box in which there is a connection between $x$ and $y$ due to asymptotic-essential connectedness. More precisely, note that if  $R(Q_{2n}(x)) < n/2$, then 
\begin{align*}
    \ell_S(x) \le |S\cap Q_{2n}(x)|+r/2,
\end{align*}
for any $n>2K$, and we can define
\begin{align*}
    N(x) := \min_{n\in \N, n>K} R(Q_{2n+1}(x))<n/2.
\end{align*}
In particular, for all $y \in Q_1(z)$ we have $\ell_S(y) \le |S\cap Q_{2N(z)+1}(z)|+r/2$ and hence
\begin{align*}
\E&\Big[\int_{S\cap Q_1} \d x\int_{S }\d y \one\{|x-y| \le \ell_S(x) + \ell_S(y)\}\Big]\\
&\le\E\Big[\sum_{n > K} \one\{N(o) = n\} |S\cap Q_1| \int_S \d y \one \{|y| \le \ell_S(y) + 1/2+r/2 + |S\cap Q_{2n+1}|\} \Big]\\
&\le\E\Big[\sum_{n > K}\sum_{m > K}\sum_{z\in \Z^2} \one\{N(o) = n\}\one\{N(z) = m\} |S\cap Q_1| \int_{S\cap Q_1(z)} \d y \one \{|y| \le \ell(y) + 1/2+r/2 + |S\cap Q_{2n+1}|\} \Big]\\
&\le\sum_{n > K}\sum_{m > K}\sum_{z\in \Z^2}\E\Big[ \one\{ N(o) = n\}\one\{N(z) = m\} |S\cap Q_1| |S\cap Q_1(z)| \one \{|z| \le |S\cap Q_{2m+1}(z)| + 1+r + |S\cap Q_{2n+1}|\} \Big],
\end{align*}
where we used $|y|\le|z|+ 1/2$ for $y \in Q_1(z)$ in order to eliminate the dependence of the integrals on $x$ and $y$ respectively and thus derive a uniform bound for the integrals.
It remains to show that this bound is integrable.
By H\"older's inequality we can bound each summand by
\begin{align*}
    \P\big(N(o) = n\big)^{1/5}&\P\big(N(o) = m\big)^{1/5} \E\big[|S\cap Q_1|^5\big]^{2/5} \P\big(|z| \le |S\cap Q_{2m+1}(z)| + 1+r + |S\cap Q_{2n+1}|\big)^{1/5}\\
    &\le C_1\exp\big(-c_1(n+m)\big)\P\big(|z| \le |S\cap Q_{2m+1}(z)| + 1+r + |S\cap Q_{2n+1}|\big)^{1/5}
\end{align*}
for some constants $C_1,c_1>0$ due to the existence of exponential moments for $|S\cap Q_1|$ as well as the assumption of exponential stabilization. 
Now, the last factor can be further bounded by Markov's inequality as
\begin{align*}
    \P\big(|z| \le   |S\cap Q_{2m+1}(z)| + 1+r + |S\cap Q_{2n+1}|\big)^{1/5}
    &\le |z|^{-3}\E\big[\big(|S\cap Q_{2m+1}(z)| + 1+r + |S\cap Q_{2n+1}|\big)^{15}\big]^{1/5}
    \\
    &\le |z|^{-3}\big(\E\big[\big(3|S\cap Q_{2m+1}|\big)^{15}] + (3+3r)^{15} + \E[(3|S\cap Q_{2n+1}|\big)^{15}]\big)^{1/5}
    \\
    &\le |z|^{-3}\big((3(2m+1)^2+3(2n+1)^2)^{15}\E\big[\big(|S\cap Q_{1}|\big)^{15}] + (3+3r)^{15}\big)^{1/5}
    \\
    &\le C_2 (mn)^{c_2}|z|^{-3},
\end{align*}
for some constants $C_2,c_2>0$. Finally, we can use $|z|>|z|_{\infty}$, where $|\cdot|_\infty$ denotes the $\ell_\infty$-norm, combined with the fact that for $k\in \N_0$ we have $\sum_{z\in Z}\one\{|z|_{\infty}=k\} \le 4k+1$, to see that
\begin{align*}
\sum_{n > K}\sum_{m > K}\sum_{k\in \N_0}(4k+1)C_2(mn)^{c_2}C_1\exp\big(-c_1(n+m)\big)k^{-3}<\infty,
\end{align*}
and hence the proof is finished. 
\end{proof}
\subsection{Proof of Theorem~\ref{Thm_Survival} }\label{ProofThm_Survival}
\begin{proof}[Proof of Theorem~\ref{Thm_Survival}]
For simplicity we only prove the case where $\muv=\delta_v$. The general case can be verified using a thinning argument where we ignore all devices that have a speed greater than $v_{\rm min}+\eps$ for an arbitrarily small $\eps$, only leading to a potentially larger value of $\lac$.

First note that the assumption that $\k$ is translation covariant and $c$-well behaved for some $c>2\varrhoImin v$ ensures that it is possible for the kernel to generate a shortest path from one street to an adjacent street while spending a sufficient amount of time on both streets to allow an infection to be past onto another device. Further, the requirement that $r>\varrhoImin v$ ensures that the communication radius is large enough such that an infection transmission can happen between devices with initial and target locations in the vicinity of the crossing. Also recall that we require that the streets of length larger than $2v\varrhoImin$ percolate with positive probability.

Since we assume $\laW$ to be arbitrary, in spirit, we think of white knights to be everywhere. We consider open streets as well as open crossings. Let us start with the streets. Let $0<\eps<\min(r/4 ,v(\varrhoImin -\varrhoWmin)/2)$ and define $n(s)= \lfloor |s|/(2r)\rfloor$, which will represent the number of devices on the street $s$ on which the infection will be propagated. We call a street $s=(c_1,c_2)$ {\em open}, if there are devices $X_0^s, \dots, X_{n(s)}^s$ such that
\begin{itemize}
\item[(S1)]$|s|>v\varrhoImin+\eps$,
\item[(S2)] for all $t$, $X_{i,t}^s \in B_\eps(c_1+\frac{ri}{2}\frac{c_2-c_1}{|c_2-c_1|})$ for all $i \in \{1,\dots, n(s)\}$,
\item[(S3)] $\rhoI(X_i^s,X_{i+1}^s)<\varrhoWmin$ for all $i \in \{1,\dots, n(s)-1\}$, and
\item[(S4)] $\rhoI(X_i^s,X_{i-1}^s)<\varrhoWmin$ for all $i \in \{2,\dots, n(s)\}$.
\end{itemize}
In words a street is open if, by Condition (S1), the street is long enough such that an infection transmission is possible. By Condition~(S2), there is a chain of devices reaching from the crossing $c_1$ to the crossing $c_2$ and such that each pair of neighboring devices has distance less than $r$. In particular, if device $X_0^s$ is infected, Condition (S3) ensures that the infection can not be stopped by white knights as they are too slow and hence the infection propagates from $X_0^s$ to $X^s_{n(s)}$.
Similarly, Condition (S4) ensures that an infection also propagates from $X^s_{n(s)}$ to $X^s_0$.
Let us call $X^s_{c_1}=X^s_1$ and $X^s_{c_2}=X^s_{n(s)}$ the {\em boundary devices} of an open street $s=(c_1,c_2)$.

Next, for a street $s=(c_1,c_2)$, denote by $[a,b]_{s,c_1}$ the segment $\{x \in s \colon x = c_1 + d\frac{c_2-c_1}{|c_2-c_1|}\text{ for some } d\in [a,b]\}$ on $s$ such that $\{0\}_{s,c_1}$ is identified with the crossing $c_1$. Consider again $0<\eps<\min(r/4 ,v(\varrhoImin -\varrhoWmin)/2)$ to be small enough such that $\vmin \varrhoImin +\eps< c/2 $. We call a crossing $c$ {\em open at time $t$} if, for any ordered pair of open streets $s_1,s_2$ that share the crossing $c$, there exists a {\em bridging device}  $X^{s_1,s_2}$ such that 
\begin{itemize}
  \item[(C1)] $X_{0}^{s_1,s_2}\in A_{s_1}=[v\varrhoImin/2+\eps/2,v\varrhoImin/2+\eps]_{s_1,c}$,
  \item[(C2)] $X^{s_1,s_2}$ has its target location in the segment $B_{s_2}=[v\varrhoImin/2,v\varrhoImin/2+\eps]_{s_2,c}$,
  \item[(C3)] $X_{t}^{s_1,s_2}\in J_{s_1}=[0,\eps/2]_{s_1,c}$ and is on its way back to its starting location,
  \item[(C4)] $\rhoI(X^{s_1,s_2},X_{c}^{s_1})< \varrhoImin+\eps/(2v)$, and
  \item[(C5)] $\rhoI(X^{s_1,s_2},X_{c}^{s_2})< \varrhoImin+\eps/v$.
  \end{itemize}
These conditions roughly guarantee that, for sufficiently long streets, the bridging device $X^{s_1,s_2}$ propagates the infection from $s_1$ to $s_2$. The time the bridging device stays on each of the streets is insufficiently long to allow to be patched by white knights. This is ensured by Condition (C1) and (C2), as the distance the devices travel on each street before entering the crossing is at most $v\varrhoImin+2\eps<v\varrhoWmin$, since $\eps$ is sufficiently small.
Together, since $s_2$ is open, then, due to the configuration of the devices on the open street and the fact that $r>v\varrhoImin $, $X^{s_1,s_2}$ will pass the infection on to the boundary device $X^{s_2}_c$ of the chain of transmitting devices on $s_2$. This is ensured by Condition (C5). Finally, Conditions~(C3) and (C4) ensure that $X^{s_1,s_2}$ becomes infected if $X^{s_1}_c$, becomes infected at time $t$, as it is on its way back to its starting location and will spend at least a time of $\varrhoWmin+\eps/(2v)$  in the vicinity of $X^{s_1}_c$.

The following statement guarantees that there exists a uniform bound for the probability that a given crossing is open for all sufficiently large times. In particular, this bound converges to one as $\la$ tends to infinity. Recall that we denoted by $H$ the range of dependence of $\k$ and define $L= {H+v\varrhoImin/2+\eps}$.
\begin{lemma}\label{lem_mixing}
For almost-all $S$ we have that $(\{c \text{ is open at time } t\})_{c\in S}$ is an independent family of events, conditioned on $S$. Further, there exists $T>0$ such that for almost-all $S$ and all crossings $c\in S$ there exists a constant $C_c=C_c(S\cap B_L(c),T)$ such that
\begin{align*}
    \P(c \text{ is open at time } t\mid S)\ge 1-\exp(-\lambda C_c)\qquad \text{ for all } t\ge T.
\end{align*}
\end{lemma}
We present the proof of the lemma later in this section. 
Consider $\tau_{s_1,c}$, the time that a boundary device associated to the crossing $c$ is infected, then, for $\tau_{s_1,c}>T$ and conditioned on $S$, we can couple the event that $c$ is open at time $\tau_{s_1,c}$ 
to a Bernoulli random variable $\mathcal{B}_c$ with parameter $1-\exp(-\lambda C_c)$, and note that the $\mathcal{B}_c$ are independent conditioned on $S$ and only depend on $S$ in the region $S\cap B_L(c)$.
For $\tau_{s_1,c}= \infty$, i.e., if the boundary device never becomes infected, we assume that the crossing is open nonetheless.
We say that a crossing is {\em open} if $\mathcal{B}_c=1$.
In order to make use of this observation, we need to ensured that the infection survives up until time $T$. 
However, note that survival for long but finite times has a positive probability and can be even constructed locally by an appropriate isolation procedure for the typical device.

Now, we consider the model at time $t\ge T$ and show that the graph of open streets and open crossings percolates with positive probability via stabilization arguments. For this, let us fix an $a$ such that $\varrhoImin <a/2 <a_{\rm c}/2$
and recall that we have the random field of stabilizing radii $\{R_x\}_{x\in\R^2}$ for the street system at our disposal, in addition to the random field $\{R^a_x\}_{x\in\R^2}$ that acts as a replacement for the asymptotically essentially connectedness of the thinned graph $S^a$.
Then, we say that $z \in \Z^2$ is \textit{$n$-open} if 
\begin{itemize}
\item[(a)] $R(Q_{3n}(nz))< n$,
\item[(b)] $R^a(Q_{n}(nz))< n/2$,
\item[(c)] every street fully contained in $Q_{3n}(nz)\cap S^a$ is open, and
\item[(d)] every crossing in $Q_{3n}(nz)$ is open.
\end{itemize} 
Condition~(b) limits the size of cells of points in $Q_{n}(nz)$ to a diameter of at most $n/2$ with respect to $S^a$. Further, under the condition, there exists a unique giant component in $Q_{n}(nz)\cap S^a$, build by the boundaries of the cells, that is connected in $Q_{3n}(nz)$. If now adjacent sites $z_1$ and $z_2$ both satisfy Condition~(b), they form a joint connected component, since there exist cells that are simultaneously in $Q_{n}(nz_1)$ and in $Q_{n}(nz_2)$.
Due to Condition~(c) and (d) the infection is able to propagate through open streets and open crossings unpatched.

Now, we can use the domination-by-product-measures theorem~\cite[Theorem 0.0]{domProd} to establish Bernoulli site percolation of $n$-open sites on $\Z^2$ since (a) guarantees $7$-dependence of $S^a$ and under Condition~(a) also the Conditions~(b) and~(c) can be decided within the $7$-dependent region. By choosing $n>L$, the crossing dependent constant $C_{\rm c}$ from Lemma~\ref{lem_mixing} is also $7$-dependent like the street system.
Hence, using translation invariance, it suffices to show that
\begin{align*}
\liminf_{n\uparrow \infty}\liminf_{\la \uparrow \infty} \P( o \text{ is $n$-open} )=1.
\end{align*}
To show this, we denote by $A(n)$, $B(n)$, $C(n,\la)$, and $D(n,\la)$, the events that the Conditions~(a), (b), (c), and (d) are not satisfied for $z=o$.
Then, $\P( o \text{ is $n$-open})\ge 1-\P(A(n)\cup B(n))-\P(C(n,\la))-\P(D(n,\la))$.
By assumption, we have that $\limsup_{n\uparrow \infty}\P(A(n)\cup B(n)) = 0$ independently of $\la$. For fixed $n$, using Markov inequality and dominated convergence, we now also have $\limsup_{\la \uparrow \infty}(\P(C(n,\la)) +\P(D(n,\la))) = 0$, as $\la$ increases, where we used that the expected number of streets and crossings in finite.

Finally, due to the translation invariance of the model, existence of an infinite cluster of open sites implies existence of an infinite cluster of good streets and good crossings and thus to a positive probability that the typical device transmits the infection to infinity. 
\end{proof}

\begin{proof}[Proof of Lemma~\ref{lem_mixing}]\label{proof_lem_mixing}
First note that the conditions only dependent on the bridging devices, which have their initial position and their destination close to the crossing and move over the crossing and are therefor independent for distinct crossings. This guarantees the independence of the events. 

In order to bound the probability of openness at late times, we compute the intensity of devices that satisfy the Conditions (1), (2) and (3) at time $t$ and then show that, as $t$ tends to infinity, this intensity converges to a positive crossing-dependent value. In turn, the probability that no such devices exist can be seen as a Poisson void probability.

Let us start by deriving expressions for the intensity of devices at time $t$ in the desired interval, i.e., devices that start on $s_1$ in the interval $A_{s_1}$ and have a target location in $B_{s_2}$, on $s_2$. Recall that the streets $s_1$ and $s_2$ have a joint crossing $c$. 
Note that, by the displacement theorem for Poisson point processes, the number of devices that satisfy the conditions is a Poisson distributed random variable and its parameter is given by $\la p_S(s_1,s_2)$, where 
\begin{align*}
p_S(s_1,s_2)=\int_{A_{s_1}}\d x \int_{B_{s_2}}\k^S(x,\d z)
\one\{0\le r_t(|x,z|_S)-|z,c|_S \le \eps/2\}, 
\end{align*}
where $|z,x|_S$ denotes the shortest distance between $x$ and $z$ on $S$ and $r_t(|z,x|_S)$ is uniquely defined as $vt=|z,x|_S n_t+r_t(|z,x|_S)$ with $n_t$ the largest odd number such that $0\le r_t(|z,x|_S)\le 2|z,x|_S$. In words, $n_t(|z,x|_S)$ denotes the number of times that the device has traveled from $x$ to $z$ and back, and since we want the device to be on the way back, we require $n_t(|z,x|_S)$ to be odd. Then, $r_t(|z,x|_S)$ is the distance from $z$ on the way back. Subtracting the distance from $z$ to the joint crossing $c$, we obtain the traveled distance on $s_1$, which has to lie in the interval $(0,\eps)_{s_1,c}$.

We will show that there is a time $T$ such that for all $t \ge T$ we can bound $p_S(s_1,s_2)$ from below by a constant that does only depend on $S\cap B_{L}(c)$. 

Let us define for $x\in A_{s_1}$ by $0<\eta_x=\eta_x(s_1,s_2,S\cap B_{L}(x),\k^S(x,\d y))= \inf_{y \in B_{s_2}}\k^S (x,y)$ the minimal density for the kernel, starting in $x$ and having a destination in $B_{s_2}$. This value $\eta_x$ is positive as the kernel is $c$-continuous, $B_{s_2} \subset B_c(x)$, and we use our assumption that $2v \varrhoImin +2\eps< c$.

In order to simplify the notation for $A_{s_1}$ and $B_{s_2}$ let us abbreviate $b=v\varrhoWmin/2$ and associate bijectively the interval $A=[b+\eps/2,b+\eps]$ to $A_{s_1}$ such that $x\mapsto {x}_{s_1,c}$. We associate $B=[b,b+\eps]$ in a similar way to $B_{s_2}$.
 Then, 
\begin{align*}
p_S(s_1,s_2)\ge \int_{A}\d x \eta_x\int_{B}\d z
\one\{0\le r_t(x+z)-z \le \eps/2\}
\end{align*}
where we also dropped the singular part in $\k^S(x,\d y)$.
Now, let us define by
\begin{align*}
\phi_x(s)=\lim_{a \downarrow 0}\frac{1}{a}\int_{B}\d z\one\{0\le r_s(x+z)-z\le a\} 
\end{align*}
the density of points points started at $x$ that are on their way back from their target location, and that are at time $s$ precisely at the origin. Note, that a device satisfies Condition~(3) at time $t$ if and only if it passes the crossing on its way to its starting location in the time interval $[t-\eps(2v)^{-1},t]$. Therefore,  we can rewrite
\begin{align}
\int_{A}\d x \eta_x\int_{B}\d z
\one\{0\le r_t(x+z)-z \le \eps/2\}=\int_{t-\eps/(2v)}^t\d s\int_{A}\d x \eta_x\phi_x(s), \label{eq_mixing_lemma_3}
\end{align}
and we will prove that $\P$-almost surely, conditioned on $S$, there exist constants $T$ and $C$ that depend on $v, \varrhoImin$, and  $\eps$ but neither on $S$, $x$ nor $\kappa^S$, such that $\inf_{t\ge T}\phi_x(t)>C>0$. 
In order to derive this statement, the main tool is the following identity for $\phi_x(t)$
\begin{align}\label{Eq_phi}
\phi_{x}(t) = \sum_{n=1}^\infty (2\eps n)^{-1} \one\{x+2b+(n-1)2(b+x)\le vt<x+2b+2\eps+(n-1)2(x+b+\eps)\}.
\end{align}
Let us explain this formula. 
Here, $n$ represents the number of of times that a device, started in $x$ reaches its destination, which is uniformly distributed in $B$.
With every reflection in $B$, the density is distributed to an interval of length $2\eps n$, where $\eps$ is the length of $B$, again uniformly at random. 
Now, $x+2b$, respectively $x+2b+2\eps$ are the shortest, respectively the longest, distance that the device, started in $x$, travels before it re-visits the crossing $c$ for the first time. For each consecutive visit in the origin an additional length of $2(b+x)$, respectively $2(x+b+\eps)$ has to be added as the length of a complete cycle from the origin to its destination and back.

In other words, each indicator becomes non-zero in $x+2b+(n-1)2(x+b)$ and adds a mass of $(2\eps n)^{-1}$ over a length of $2n\eps$. For $n<(x+b)/\eps$ the indicator functions are disjoined, i.e. for $vt<x+2b+((x+b)/\eps-1)2(x+b)$ there is at most one indicator function that is non-zero.
Now, for $n>(x+b)/\eps$, the interval for $vt$ in which the indicator function is non-negative is longer than the minimal cycle length. In other words, it is possible that at least two indicator functions are non-negative.

Now, for $(x+b)/\eps<n<2(x+b)/\eps$ the indicator functions contribute for more than the length of a full cycle $2(x+b)$ and less than twice that length.
Therefore, for each $vt \in [x+2b+(n_1-1)2(b+x),x+2b+(n_2-1)2(b+x)]$ there are either exactly one or two non-negative indicators, where $n_i=\lceil i(x+b)/\eps \rceil$. As each indicator has a weight of at least $(4\epsilon(x+b)/\eps)^{-1}$ we can bound $\phi_x(t)>(4\epsilon(x+b)/\eps)^{-1}$.

More general, in the interval $vt \in [x+2b+(n_m-1)2(b+x),x+2b+(n_{m+1}-1)2(b+x)]$
there are either $m$ or $m+1$ overlapping indicator functions whose contribution can be each be bounded by $(2n_{m+1}\eps)^{-1}$. Therefore, we can bound  
\begin{align}\label{eq_lem_mixing_2}
    \phi_x(t)>m(2(n_{m+1})\eps)^{-1} = \frac{m}{2\eps \lceil (m+1) (x+b)/\eps \rceil} \ge \frac{1}{2\eps \lceil 2(x+b)/\eps \rceil}\ge \frac{1}{4(2b+\eps)}=C>0,
\end{align}
where we used that the bound monotonously increases in $m$ and $x\le b+\eps$.
Therefore, for $T=(x+2b+(n_1-1)2(b+x))v^{-1}$ the desired bound is realized and integration of~\ref{eq_mixing_lemma_3} yields
\begin{align*}
    p_S(s_1,s_2)\ge \eps(2v)^{-1}(8b+4\eps)^{-1} \int_{A}\d x \eta_x=\tilde p_S(s_1,s_2)>0. 
\end{align*}

Noting that the requirements in Conditions~(4) and (5) can be considered independently, respectively with probability $p_4=\varrhoI([\varrhoImin,\varrhoImin+\eps/2])$ and $p_5=\varrhoI([\varrhoImin,\varrhoImin+\eps])$, the number of devices that satisfy the criterion of a bridging device $X^{s_1,s_2}$ can be bounded from below by an Poisson-distributed random variable with expectation $\tilde p_S(s_1,s_2)p_4p_5$.
As $c$ is open if there is at least one bridging device for each combination of open adjacent streets, we can bound
\begin{align*}
    \P(c \text{ is open at time } t\mid S)\ge \prod_{\substack{i\neq j, s_i,s_j \text{ open}\\c \in s_i,c \in s_j}} \big(1-\exp(-\lambda \tilde p(s_i,s_j)p_4p_5)\big)\qquad \text{ for all } t\ge T,
\end{align*}
which proves the lemma.
\end{proof}
\subsection{Proof of Theorem~\ref{Thm_Extinction} }\label{Proof_Thm_Extinction}
For the proof we adopt a multiscale argument in the spirit of the proof of existence of subcritical regimes for the Poisson--Boolean model with random radii in~\cite{G08}. We note that this approach has been also applied in the setting of Cox--Boolean models with random radii in~\cite{jahnel2022phase}. More specifically we adapt the arguments for absence of percolation in the connectivity graph for sufficiently large velocities presented in~\cite{HJCW22}. Roughly speaking, we prove that a Cox--Boolean model in which any two devices, that could transmit the infection somewhere on their path and are sufficiently close, are connected, becomes subcritical if the white-knight intensity is sufficiently large. 

First, we call a street $s=(c_1,c_2)$ {\em blocked for (device) $x$} if either $|s|<v_{\rm min}\varrhoImin/2$ or if there exist two white knights $Y_1,Y_2\in Y^{\laW}$ such that for $i\in \{1,2\}$ we have that $\varrhoImin>\rhoW(Y_i,x)$, $Y_{i,t}\in s$ and $|c_i-Y_{i,t}|<r$ for all $t\ge 0$.
Otherwise we call $s$ {\em unblocked for $x$}.
Then, we consider the thinned process 
$$X^{\la,{\rm th}}=\{X_i\in X^\la\colon \text{there exists }y\in\supp(\k^S(X_i,\d y))\text{ and } s\in\ell_S(X_i,y)\text{ such that $s$ is unblocked for $X_i$}\}.$$
In words, for devices in $X^\la\setminus X^{\la,{\rm th}}$ every possible path only contains streets that are either to short to perform any transmission of the infection or that contain a white knight that patches the device before it can transmit the infection on the street. Let us note that here we use that $r>v_{\rm max}\varrhoWmin$ is large enough to allow a successful patch.  Using these definitions, we note that, if the infection survives globally, then this implies that the geostatistical Cox--Boolean model
$$\C=\bigcup_{X_i\in X^{\la,{\rm th}}}B_{\ell_S(X_i)}(X_i),$$
contains an unbounded component. Hence, in order to prove absence of global survival, it suffices to prove absence of percolation in this Cox--Boolean model, which is also well defined by the local-connectedness assumption.
Let $\C_x(V)$ denote the connected component containing $x$ of the geostatistical Cox--Boolean model, based on points in $V\subset\R^2$. Then, we define for any $x\in \R^2$ and $\alpha>0$ the event
\begin{align*}
G(x,\a)&=\big\{\C_x(B_{10\a}(x))\not\subset B_{8\a}(x)\big\}
\end{align*}
that the cluster of $x$, only using points in $B_{10\a}(x)$, reaches beyond $B_{8\a}(x)$. Consider the set of $\a$-local points
\begin{align*}
A_S(\a)=\{x\in S\colon \ell_S(x)<\a \},
\end{align*}
let $S^*$ denote the Palm version of $S$ and recall the stabilization radii $R_y$ of the street system. Recall that $R(V)=\sup_{y\in V\cap\, \Q^2}R_y$ for any $V\subset\R^2$. Then, we can employ the following key lemma that establishes a scaling relation in the model.
\begin{lemma}[{\cite[Lemma III.1]{HJCW22}}]\label{Lem_A}
Consider the geostatistical Cox--Boolean model. Then, there exists a constant $c>0$ such that for all $\a>0$, we have 
\begin{align*}
\P\big(G(o,10\a)\big)&\le c\P \big(G(o,\a)\big)^2+c\a^2\P\big(o\in A^{\rm c}_{S^*}(\a)\big)+c\P\big(R(Q_{10\a})\ge \a\big)
\end{align*}
\end{lemma}
In the next lemma we deviate from~\cite[Proof of Theorem II.3]{HJCW22}. We show that the local percolation probability becomes zero for large white-knight intensities. 
\begin{lemma}\label{Lem_C}
There exists $c>0$ such that 
$$\P\big(G(o,\a)\big)\le c \a^2\P\big(\text{there exists }y\in\supp(\k^{S^*}(o,\d y))\text{ and at least one unblocked }s\in\ell_{S^*}(o,y)\big)$$ 
and $\P\big(\text{there exists }y\in\supp(\k^{S^*}(o,\d y))\text{ and at least one unblocked }s\in\ell_{S^*}(o,y)\big)$ tends to zero as $\laW$ tends to infinity. 
\end{lemma}
\begin{proof}[Proof of Lemma~\ref{Lem_C}]
Note that 
\begin{align*}
\P\big(G(o,\a)\big)&\le \P\big(X^{\la,{\rm th}}(B_{10\a})>0\big)\le \E\big[X^{\la,{\rm th}}(B_{10\a})\big]\\
&\le \la\E\Big[\int_{S\cap B_{10\a}}\d x\P\big(\text{there exists }y\in\supp(\k^S(x,\d y))\text{ and } s\in\ell_S(x,y)\text{ such that $s$ is unblocked for $x$}\big)\Big]\\
&\le \la 10^2\a^2\P\big(\text{there exists }y\in\supp(\k^{S^*}(o,\d y))\text{ and } s\in\ell_{S^*}(o,y)\text{ such that $s$ is unblocked for $o$}\big)
\end{align*}
But, by dominated convergence $\P\big(\text{there exists }y\in\supp(\k^{S^*}(o,\d y))\text{ and } s\in\ell_{S^*}(o,y)\text{ such that $s$ is unblocked for $o$}\big)$ tends to zero as $\laW$ tends to infinity and hence we arrive at the desired result. 
\end{proof}
As a final input for the proof, we recall the following essential result about convergence properties of functions satisfying some scaling inequality. 
\begin{lemma}[{\cite[Lemma 3.7]{G08}}]\label{Lem_F}
Let $f$ and $g$ be two bounded measurable functions from $[1,\infty]$ to $[0,\infty)$. Additionally, let $f$ be bounded by $1/2$ on $[1,10]$, $g$ be bounded by $1/4$ on $[1,\infty]$ and assume 
\begin{align*} 
f (\a) \le f (\a/10)^2 + g(\a),\qquad\text{for all }\a\ge 10.
\end{align*}
Then, $\lim_{\a\uparrow\infty}g(\a)=0$ implies that $\lim_{\a\uparrow\infty}f(\a)=0$. 
\end{lemma}
Now we are in the position to prove the theorem.
\begin{proof}[Proof of Theorem~\ref{Thm_Extinction}]
Recall that $S$ is assumed to be a PVT or PDT. In order to prove that $\laWc<\infty$, it suffices to show that $\lim_{\a\uparrow\infty}\P(X^{\la,{\rm th}}(\C_o(\R^2))> X^{\la,{\rm th}}(B_{8\a}))=0$ for all sufficiently large $\laW$. But this is true if $\lim_{\a\uparrow\infty}\P\big(G(0,\a)\big)=0$, since by~\cite{HJCW22}[Lemma III.4 and III.5], we have that 
\begin{align*}
\P\big(X^{\la,{\rm th}}(\C_o(\R^2))> X^{\la,{\rm th}}(B_{8\a})\big)&\le \P\big(G(o,\a)\big)+\la\E\Big[\int_{S}\d x\one\{10\a\le |x|\le 9\a+\ell_S(x)\}\Big], 
\end{align*}
where $\E\Big[\int_{S}\d x\one\{10\a\le |x|\le 9\a+\ell_S(x)\}\Big]$ tends to zero as $\a$ tends to infinity. Now, in order to show $\lim_{\a\uparrow\infty}\P\big(G(0,\a)\big)=0$ we apply Lemma~\ref{Lem_A}, Lemma~\ref{Lem_C} and Lemma~\ref{Lem_F} for proper choices of $f$ and $g$. 

For this, note that, since we assume stabilization, in Lemma~\ref{Lem_A} we have that $\P(R(Q_{10\a})\ge \a)$ tends to zero as $\a$ tends to infinity and the same is true for $ \a ^2\P\big(o\in A^{\rm c}_{S^*}(\a)\big)$ by~\cite{HJCW22}[Lemma III.2]. Hence, we can define 
\begin{align*}
\a_1&:=\inf\{s\ge 1\colon \phi_1(\a)=\a ^2\P\big(o\in A^{\rm c}_{S^*}\big)<(8c^2)^{-1}\text{ for all }\a\ge s\}\\
\a_2&:=\inf\{s\ge 1\colon \phi_2(\a)=\P\big(R(Q_{10\a})\ge \a\big)<(8c^2)^{-1}\text{ for all }\a\ge s\},
\end{align*}
and set $\a_{\rm c}=\a_1\vee\a_2$ and note that $\a_{\rm c}<\infty$. Further, we let
\begin{align*}
\laWc&=\inf\{\laW\ge 1\colon \P\big(\text{there exists }y\in\supp(\k^{S^*}(o,\d y))\text{ and at least one unblocked }s\in\ell_{S^*}(o,y)\big)<\frac 1 2 (100c\a_{\rm c})^{-2}\},
\end{align*} 
where $\laWc<\infty$ by Lemma~\ref{Lem_C}. Then, we denote
\begin{align*}
f(\a)=c\P\big(G(o,10\a_{\rm c}\a)\big)\quad&\text{ and }\quad g(\a)=c^2\big(\phi_1(\a_{\rm c}\a)+\phi_2(\a_{\rm c}\a)\big)
\end{align*}
and hence, using again Lemma~\ref{Lem_C}, we have that 
\begin{align*}
f(\a)&\le 1/2,\qquad\text{ for all }1\le \a\le10\text{ and }\laW>\laWc,
\end{align*}
and by definition,
\begin{align*}
g(\a)&\le 1/4,\qquad\text{ for all }1\le \a.
\end{align*}
Finally, using Lemma~\ref{Lem_A}, we have that 
\begin{align*}
f(\a)&\le c^2\big(\P\big(G(o,\a_{\rm c}\a)\big)^2+\phi_1(\a_{\rm c}\a)+\phi_2(\a_{\rm c}\a) \big)=f(\a/10)^2+g(\a)\quad\text{ for all }\a \ge 10.
\end{align*}
Hence, since $\lim_{\a\uparrow\infty}g(\a)=0$, an application of Lemma~\ref{Lem_F} gives the result. 
\end{proof}

\subsection{Proof of Theorem~\ref{thm_InOut} }\label{Proof_thm_InOut}
We prove the theorem in two parts. 
\begin{proof}[Proof of Theorem~\ref{thm_InOut} Part (2) and (3)]
Note that Part (2) is the statement of Theorem~\ref{Thm_Survival} which puts no restrictions in the white-knight intensity. Further, Part (3) follows from \cite{HJCW22}[Theorem 3.3] where  a sub-critical percolation regime is established for sufficiently fast speeds, independent of the device intensities. This in particular implies that global survival is possible only with probability zero. 
\end{proof}
Hence, it only remains to prove Part (1). The main idea is that small velocities make it impossible to transmit the infection through crossings if sufficiently many white knights are available. We say that a street $s=(c_1,c_2)$ is {\em closed} if 
\begin{itemize}
\item[(A)] we have that $|s|>2v\rhoW$,
\item[(B)] there exist white knights $Y^s\subset Y^{\laW}$ such that for all $t\ge 0$ and $x\in s$ we have $|x-Y_{i,t}|<r/2$ for some $Y_i\in Y^s$, and
\item[(C)] if at time $t\ge 0$, there exists an infected $X_i\in X^\la$ with $X_{i,t}=c_1$ (respectively $X_{i,t}=c_2$), then, $\{X_j\in X^\la\cap s\colon |X_{j,u}-c_2|\le v\rhoW,\text{ for some }u\in \tau(t)\}=\emptyset$ (respectively $\{X_j\in X^\la\cap s\colon |X_{j,u}-c_1|\le v\rhoW,\text{ for some }u\in \tau(t)\}=\emptyset$). Here, $\tau(t)$ is the time interval in which the infection reaches $s\cap B_{v\rhoW+r}(c_2)$ (respectively in $s\cap B_{v\rhoW+r}(c_1)$).
\end{itemize}
A street is called {\em open} if it is not closed.
The key idea is that a closed street cannot be used for the transmission of the infection. Indeed, Condition~(A) guarantees that the street is sufficiently long in order to allow the white knights to patch, and the probability that the condition is satisfied tends to one as $v$ tends to zero. Further, Condition~(B) guarantees the existence of white knights at all times close to every point on the street and in particular close to the crossings. These white knights help to patch any incoming infected device and the condition becomes likely to be satisfied for large white-knight intensities. Now, for Condition~(C), note first that a large device intensity $\la$ makes it probable that an infection travels through a street unpatched since we assume that $\rhoW>\rhoI$. Also, due to the continuous movement of the devices, it is unavoidable that an infected device enters the street at some time. However, in order for an infected device to leave the street unpatched, it needs to receive the infection at the opposite end of the street, sufficiently close to the crossing, i.e., in $s\cap B_{v\rhoW+r}(c_2)$ (respectively $s\cap B_{v\rhoW+r}(c_1)$). Only in this case it can exit the street unpatched from the blocking white knight that is positioned there. Very roughly then, again for small $v$, Condition~(C) becomes likely. 

What makes the probability that Condition~(C) is satisfied challenging to control is the fact that, in principle, a street at any given time could host devices from distant spatial areas. In order to cope with this we again employ a multi-scale argument very similar to the one presented in the previous Section~\ref{Proof_Thm_Extinction}. 
Using the definitions for $\ell_S$ from Section~\ref{Proof_Thm_Extinction}, we now consider the thinned point cloud
$$X^{\la,{\rm th}}=\{X_i\in X^\la\colon \text{there exists }y\in\supp(\k^S(X_i,\d y))\text{ and } s\in\ell_S(X_i,y)\text{ such that $s$ is open}\}.$$
Again, the idea is that any device that is not in $X^{\la,{\rm th}}$ only finds closed streets on any possible path and hence can not contribute to the spread of the infection. Thus, it suffices to show subcriticality in the associated geostatistical Cox--Boolean model
$$\C=\bigcup_{X_i\in X^{\la,{\rm th}}}B_{\ell_S(X_i)}(X_i).$$
However, the situation is more complicated compared to the proof of Theorem~\ref{Thm_Extinction} since the thinning procedure depends on the entire point cloud. Still, we can employ the Lemmas~\ref{Lem_A} and \ref{Lem_F} and only need to replace the Lemma~\ref{Lem_C} by the following statement, which is presented in terms of $\P^*$, the Palm version of $\P$.
\begin{lemma}\label{Lem_C2}
There exists $c>0$ such that 
$$\P\big(G(o,\a)\big)\le c \a^2\P^*\big(\text{there exists }y\in\supp(\k^S(o,\d y))\text{ and } s\in\ell_S(o,y)\text{ such that $s$ is open}\big)$$ 
and 
$$\limsup_{\laW\uparrow\infty}\limsup_{v\downarrow0}\P^*\big(\text{there exists }y\in\supp(\k^S(o,\d y))\text{ and } s\in\ell_S(o,y)\text{ such that $s$ is open}\big)=0.$$ 
\end{lemma}
\begin{proof}[Proof of Lemma~\ref{Lem_C2}]
First note that we can use the Slivnyak--Mecke formula and translation invariance to estimate 
\begin{align*}
\P\big(G(o,\a)\big)&\le \P\big(X^{\la,{\rm th}}(B_{10\a})>0\big)\le \E\big[X^{\la,{\rm th}}(B_{10\a})\big]\\
&\le \la 10^2\a^2\P^*\big(\text{there exists }y\in\supp(\k^{S}(o,\d y))\text{ and } s\in\ell_{S}(o,y)\text{ such that $s$ is open}\big).
\end{align*}
We highlight here that, different to the situation of Lemma~\ref{Lem_C}, the probability also extends to the devices $X^\la$ since openness is defined not only with respect to white knights and the street system. 

In order to show the second statement of the lemma, we employ stabilization based on the stabilizing random field $\{R_x\}_{x\in\R^2}$. For all sufficiently large $n$, we have that for $\eps'>0$,
\begin{align*}
\P^*\big(&\text{there exists }y\in\supp(\k^{S}(o,\d y))\text{ and } s\in\ell_{S}(o,y)\text{ such that $s$ is open}\big)\\
&\le \P^*\big(R(Q_{2n})<n/2, \text{there exists }y\in\supp(\k^{S}(o,\d y))\text{ and } s\in\ell_{S}(o,y)\text{ such that $s$ is open}\big)+\eps'\\
&\le \P^*\big(\text{there exists an open $s\in Q_{2n}$}\big)+\eps',
\end{align*}
where we used that $\P(R(Q_{2n})\ge n/2)$ tends to zero as $n$ tends to infinity, as well as the fact that under the event $R(Q_{2n})< n/2$ any shortest path starting in $o$ must be contained in $Q_{2n}$. Now, by the Markov inequality, 
\begin{align*}
\P^*\big(\text{there exists an open $s\in Q_{2n}$}\big)\le \E^*\Big[\sum_{s\in S\cap Q_{2n}}\P^*(s\text{ is open}|S)\Big]
\end{align*}
and, using dominated convergence and the fact that $\E^*[\#\{s\in S\cap Q_{2n}\}]<\infty$ for PVT and PDT, it suffices to show that, for all $s\in S\cap Q_{2n}$, 
$$\limsup_{\laW\uparrow\infty}\limsup_{v\downarrow0}\P^*(s\text{ is open}|S)=0,$$
almost-surely with respect to $S$.
In view of the definition of closed streets, let $C(v)$ denote the event that Condition~(A) is satisfied for $s$ and similar $D(\laW)$ and $E(v,\la,\laW)$ for the events associated to the Conditions~(B) and (C). Then, 
\begin{align*}
\P^*(s\text{ is open}|S)\le \P(C^c(v)|S)+\P( D^c(\laW)|S)+\P( E^c(v,\la,\laW)|S), 
\end{align*}
where the first and third summands tend to zero as $v$ tends to zero and then also the second summand tends to zero as $\laW$ tends to infinity. This finishes the proof.
\end{proof}

\begin{proof}[Proof of Theorem~\ref{thm_InOut} Part 1]
For the proof we can follow the same line of arguments as used for the proof of Theorem~\ref{Thm_Extinction} above, only replacing Lemma~\ref{Lem_C} by Lemma~\ref{Lem_C2} and subsequently adjusting the arguments. We need to prove existence of $\laWc<\infty$ such that for all $\laW>\laWc$ there exists $v_{c}(\la,\laW)$ such that $\lim_{\a\uparrow\infty}\P(X^{\la,{\rm th}}(\C_o(\R^2))> X^{\la,{\rm th}}(B_{8\a}))=0$  for all $v<v_{c}(\la,\laW)$. Referencing to~\cite{HJCW22}[Lemma III.2, Lemma III.4 and III.5] as well as to Lemma~\ref{Lem_A} and the initial argument presented in Theorem~\ref{Thm_Extinction} it suffices to note that, with the help of Lemma~\ref{Lem_C2}, there exists $\laWc<\infty$ such that for all $\laW>\laWc$ there exists $v_{c}(\la,\laW)$ such that for all $v<v_{c}(\la,\laW)$ 
$$
\P^*\big(\text{there exists }y\in\supp(\k^S(o,\d y))\text{ and } s\in\ell_S(o,y)\text{ such that $s$ is open}\big)<\frac 1 2 (100c\a_{\rm c})^{-2},
$$
where $\a_{\rm c}$ is defined as in the proof of Theorem~\ref{Thm_Extinction}. Then, using the same definitions for $f$ and $g$, the result follows by an application of Lemma~\ref{Lem_F}. 
\end{proof}
\subsection{Proof of Proposition~\ref{Prop_FPP_fixed_T}}\label{sec_Prop4-1}
\begin{proof}[Proof of Proposition~\ref{Prop_FPP_fixed_T}]
First note that by the assumption of supercriticality there is a positive probability that the typical device is connected to infinity within $g_{\rho,r}(S,X^\la)$ via at least one path of susceptible devices, i.e., in the absence of white knights. However, in case that $T_{\rm I}\le T_{\rm W}$, any white knight that is connected to the path can only eliminate the infection after the infected device has already passed the infection towards the subsequent susceptible device, leading to an unstoppable sequence of infections along the path, independent of the choice of $\laW$.

On the other hand, in case $T_{\rm I}>T_{\rm W}$, it suffices to consider the model conditioned on the event $\{X_o\leftrightsquigarrow\infty\}$ that $X_o$ is in the infinite cluster of $g_{\rho,r}(S,X^\la)$. 
Let us denote by $D$ the graph distance of $X_o$ to the set of white knight and note that $D<\infty$ almost surely as $\laW>0$. 
Then, closest white knight $Y_i$ experiences an infection attempt at time $DT_{\rm I}$. Note that, at this time, any infected device has a distance to $Y_i$ given by at most $2D$. This comes from the fact that the transmission times are deterministic and there is a one-to-one correspondence between graph distance to the typical device and elapsed time. However, there exists $N\in \N$ such that $NT_{\rm I}>(N+2D)T_{\rm W}$ and therefore the infection can reach at most a graph distance of $N+D$. Indeed, any device $X_i$ at graph distance $N+D$ from $X_o$, is in contact with $Y_i$ via $N+2D$ edges connecting infected or patched devices. But, since the patching is faster than the infection, for this $N$, on the joint path between $X_o$, $Y_i$ and $X_i$, the patch will have caught up with the infection before time $(N+D)T_{\rm I}$ which is strictly smaller than $NT_{\rm I}$.
\end{proof}

\subsection{Proof of Theorem~\ref{Thm_Exp_Sur} }\label{proofThm_Exp_Sur}

\begin{proof}[Proof of Theorem~\ref{Thm_Exp_Sur}]
As $\varrhoWmin > 0$, we can choose $b$ large enough such that $\varrhoIminb<\varrhoWmin$.
Furthermore, we can choose $\mu'$ with $\vminp$ such that
$\vmin\rho<\vminp\varrhoIminb<\min(a_{\rm c}/2,r,c/2)$. 
Now $\vminp$ and $\varrhoIminb$ satisfy the Conditions of Theorem~\ref{Thm_Survival}.
A close inspection of the proof of 
Theorem~\ref{Thm_Survival} reveals that, under the conditions in this theorem, there exists an infinite sequence of crossing devices, boundary devices and devices connecting boundary devices $\{X_{n_i}\}_{i\in \N}$, connecting the typical device to infinity, in the dynamical model considered in Theorem~\ref{Thm_Survival}.
Let us argue that this given sequence is also connected in $g_{\rho,r}(S, X^\la)$. Indeed, as the chain of devices between boundary devices is strongly localized, they are connected for all $t\ge 0$, especially for a continuous time-window of length $\rho$.
Now, the crossing devices are localized in such a way, that they travel at least a distance of $\vminp\varrhoIminb$ on their respective streets. As $\vmin\rho<\vminp\varrhoIminb$ the bridging devices can connect to their respective boundary devices. Therefore, the infection can propagate in $g_{\rho,r}(S, X^\la)$ along the given sequence to infinity, which finishes the proof.
\end{proof}

\subsection{Proof of Theorem~\ref{Thm_Exp_Sur_2} }\label{proofThm_Exp_Sur_2}
\begin{proof}[Proof of Theorem~\ref{Thm_Exp_Sur_2}]
Using similar arguments as in the analysis of chase-escape dynamics on the Gilbert graph in~\cite{hinsen2020phase}, we call a device $X_i\in V^{M,p}$ {\em good} if, in case of an infection, the infection is transmitted to all neighbors of $X_i$ before any neighbor of $X_i$ is able to patch it. If there is an infinite path of {\em good} devices, every device on this path is guaranteed to become infected and hence, if there exists an infinite component of good devices with positive probability, this implies survival of the infection also with positive probability.
Now, if a device $X_i$ has at most $M$ neighbors, the probability to be good is bounded by
\begin{align}\nonumber
\P\Big(\min_{X_j \in N(X_i)} \rhoW(X_j,X_i) > \max_{X_j \in X^\la \cap N(X_i)} \rhoI(X_i,X_j)\Big) &\ge \P\Big(\min_{i \in \{1,\dots,M\}} \rho_{{\rm W},i} > \max_{i \in \{1,\dots, M\}} \rho_{{\rm I},i}^b\Big)
\\\label{bound_1}
&=\P\Big(\min_{i \in \{1,\dots,M\}} \rho_{{\rm W},i} > b^{-1}\max_{i \in \{1,\dots, M\}} \rho_{{\rm I},i}\Big),
\end{align}
independently of $X_i$, where the $\rho_{{\rm I},i}^b$, and respectively the $\rho_{{\rm W},i}$, are iid random variables with distribution $\varrhoIb$ and $\varrho_{\rm W}$.

Now, by assumption, we can choose $M$ and $p$ such that  $g^{M,p}_{\rho,r}(S,X^\la)$ percolates with positive probability.
Hence, as $\varrhoW$ had no atom at zero
we can choose $b=b(M)$ sufficiently large such that the bound~\eqref{bound_1} is larger than $p$. 

Note that the bound is uniform over the devices, since the infection and patch times are independent by definition, and therefore we can couple the connectivity graph to an independent thinning with parameter $p$ that dominates the process of good devices. 
As a consequence, $g^{M,p}_{\rho,r}(S,X^\la)$ is a subgraph of the thinned graph of good devices with at most $M$ neighbors. Hence, we have established percolation of good devices with positive probability as desired. 
\end{proof}

\subsection{Proof of Theorem~\ref{Thm_Exp_Ext} }\label{ProofThm_Exp_Ext}
Our strategy is again to employ the arguments from Section~\ref{Proof_Thm_Extinction} and Section~\ref{Proof_thm_InOut}. That is, we realize that it suffices to show absence of percolation in the geostatistical Cox--Boolean model
$$\C=\bigcup_{X_i\in X^{\la,{\rm th}}}B_{\ell_S(X_i)}(X_i),$$
where, 
$$X^{\la,{\rm th}}=\{X_i\in X^\la\colon \min_{X_j \in N_\la(X_i)} \rhoI(X_i,X_j) < \min_{X_j \in N_\laW(X_i)} \rhoW(X_j,X_i)\},$$
with $N_\la(X_i)$ the set of neighbors of $X_i$ in $g_{\rho,r}(S,X^\la)$ and $N_{\laW}(X_i)$ the set of neighbors of $X_i$ in $g_{\rho,r}(S,X_i\cup Y^\laW)$. Indeed, devices in $X^\la\setminus X^{\la,{\rm th}}$ are patched before they can transmit the infection and can therefore not contribute to the dissipation of the infection.

Again, the thinning procedure depends on the entire point cloud, but, we can still employ Lemma~\ref{Lem_A} and Lemma~\ref{Lem_F} from above and only need to replace Lemma~\ref{Lem_C} by the following statement, which is presented in terms of $\P^*$, the Palm version of $\P$.
\begin{lemma}\label{Lem_C3}
There exists $c>0$ such that 
$$\P\big(G(o,\a)\big)\le c \a^2\P^*\big(\min_{X_j \in N_\la(o)} \rhoI(o,X_j) < \min_{X_j \in N_\laW(o)} \rhoW(X_j,o)\big)$$ 
and 
$$\limsup_{\laW\uparrow\infty}\P^*\big(\min_{X_j \in N_\la(o)} \rhoI(o,X_j) < \min_{X_j \in N_\laW(o)} \rhoW(X_j,o)\big)=0.$$ 
\end{lemma}
\begin{proof}[Proof of Lemma~\ref{Lem_C3}]
Again, by the Slivnyak--Mecke formula and translation invariance we can estimate
\begin{align*}
\P\big(G(o,\a)\big)&\le \P\big(X^{\la,{\rm th}}(B_{10\a})>0\big)\le \E\big[X^{\la,{\rm th}}(B_{10\a})\big]\\
&\le \la 10^2\a^2\P^*\big(\min_{X_j \in N_\la(o)} \rhoI(o,X_j) < \min_{X_j \in N_\laW(o)} \rhoW(X_j,o)\big).
\end{align*}
Now, 
\begin{align}\label{Est_1}
\P^*\big(\min_{X_j \in N_\la(o)} \rhoI(o,X_j) < \min_{X_j \in N_\laW(o)} \rhoW(X_j,o)\big)\le \E^*\big[\P(\rhoW>\varrhoImin)^{N_\laW(o)}\big],
\end{align}
where $\rhoW$ is a random variable distributed according to $\varrhoW$. Note that $\P(\rhoW>\varrhoImin)<1$ by our assumptions on $\varrhoWmin$. Further, the expression on the right-hand side of~\eqref{Est_1} can be further bounded from above by considering only those white knights that mimic the behavior of the typical device at the origin in the sense that they start close to the origin and have a target location close to the target location of the typical device. Then, by dominated convergence, the right-hand side of~\eqref{Est_1} tends to zero as $\laW$ tends to infinity and this finishes the proof.
\end{proof}
As now the Lemma is proven this concludes the proof of Theorem~\ref{Thm_Exp_Ext}.

\begin{proof}[Proof of Theorem~\ref{thm_InOut} Part 1]
For the proof we can follow the same line of arguments as used for the proof of Theorem~\ref{Thm_Extinction} above, only replacing Lemma~\ref{Lem_C} by Lemma~\ref{Lem_C2} and subsequently adjusting the arguments. We need to prove existence of $\laWc<\infty$ such that for all $\laW>\laWc$ there exists $v_{c}(\la,\laW)$ such that $\lim_{\a\uparrow\infty}\P(X^{\la,{\rm th}}(\C_o(\R^2))> X^{\la,{\rm th}}(B_{8\a}))=0$  for all $v<v_{c}(\la,\laW)$. Referencing to~\cite{HJCW22}[Lemma III.2, Lemma III.4 and III.5] as well as to Lemma~\ref{Lem_A} and the initial argument presented in Theorem~\ref{Thm_Extinction} it suffices to note that, with the help of Lemma~\ref{Lem_C2}, there exists $\laWc<\infty$ such that for all $\laW>\laWc$ there exists $v_{c}(\la,\laW)$ such that for all $v<v_{c}(\la,\laW)$ 
$$
\P^*\big(\text{there exists }y\in\supp(\k^S(o,\d y))\text{ and } s\in\ell_S(o,y)\text{ such that $s$ is open}\big)<\frac 1 2 (100c\a_{\rm c})^{-2},
$$
where $\a_{\rm c}$ is defined as in the proof of Theorem~\ref{Thm_Extinction}. Then, using the same definitions for $f$ and $g$, the result follows by an application of Lemma~\ref{Lem_F}. 
\end{proof}

\section*{Acknowledgements}
This work was funded by Orange Labs S.A.~under the project ID CRE K07151, the German Research Foundation under Germany's Excellence Strategy MATH+: The Berlin Mathematics Research Center, EXC-2046/1 project ID 390685689, and the German Leibniz Association via the Leibniz Competition 2020.


\bibliography{wias}
\bibliographystyle{alpha}

\end{document}